\documentclass[12pt]{article}
\usepackage{amsfonts,amsmath}
\def \beq {\begin{eqnarray}}
\def \eeq {\end{eqnarray}}
\def \beqn {\begin{eqnarray*}}
\def \eeqn {\end{eqnarray*}}

\newcommand{\halmos}{\rule{1ex}{1.4ex}}

\newcounter{for}[section]

\newtheorem{itlemma}{Lemma}[section]
\newtheorem{itproposition}[itlemma]{Proposition}
\newtheorem{itconjecture}[itlemma]{Conjecture}
\newtheorem{theorem}[itlemma]{Theorem}
\newtheorem{itcorollary}[itlemma]{Corollary}
\newtheorem{itremark}[itlemma]{Remark}
\newtheorem{itremarks}[itlemma]{Remarks}
\newtheorem{itdefinition}[itlemma]{Definition}
\newtheorem{itexample}[itlemma]{Example}

\newenvironment{fact}{\begin{itfact}\rm}{\end{itfact}}
\newenvironment{claim}{\begin{itclaim}\rm}{\end{itclaim}}
\newenvironment{lemma}{\begin{itlemma}}{\end{itlemma}}
\newenvironment{remark}{\begin{itremark}}{\end{itremark}}
\newenvironment{remarks}{\begin{itremarks}\rm}{\end{itremarks}}
\newenvironment{corollary}{\begin{itcorollary}}{\end{itcorollary}}
\newenvironment{proposition}{\begin{itproposition}}{\end{itproposition}}
\newenvironment{conjecture}{\begin{itconjecture}}{\end{itconjecture}}
\newenvironment{definition}{\begin{itdefinition}\rm}{\end{itdefinition}}
\newenvironment{example}{\begin{itexample}\rm}{\end{itexample}}
\newenvironment{proof}{\noindent {\em Proof}.\ \
}{\hspace*{\fill}$\halmos$\medskip}
\newcommand{\be}[1]{\addtocounter{for}{1} \begin{equation}\label{#1}}
\newcommand{\ee}{\end{equation}}
\newcommand{\bl}[1]{\begin{lemma}\label{#1}}
\newcommand{\br}[1]{\begin{remark}\label{#1}}
\newcommand{\brs}[1]{\begin{remarks}\label{#1}}
\newcommand{\bt}[1]{\begin{theorem}\label{#1}}
\newcommand{\bd}[1]{\begin{definition}\label{#1}}
\newcommand{\bp}[1]{\begin{proposition}\label{#1}}
\newcommand{\bcon}[1]{\begin{conjecture}\label{#1}}
\newcommand{\bc}[1]{\begin{corollary}\label{#1}}
\newcommand{\bfact}[1]{\begin{fact}\label{#1}}
\newcommand{\bex}[1]{\begin{example}\label{#1}}
\newcommand{\ec}{\end{corollary}}
\newcommand{\efact}{\end{fact}}
\newcommand{\eex}{\end{example}}
\newcommand{\el}{\end{lemma}}
\newcommand{\er}{\end{remark}}
\newcommand{\ers}{\end{remarks}}
\newcommand{\et}{\end{theorem}}
\newcommand{\ed}{\end{definition}}
\newcommand{\ep}{\end{proposition}}
\newcommand{\econ}{\end{conjecture}}
\newcommand{\epr}{\end{proof}}
\newcommand{\bpr}{\begin{proof}}
\newcommand{\bcl}[1]{\begin{claim}\label{#1}}
\newcommand{\ecl}{\end{claim}}
\newcommand{\ecs}{\end{corollary}}
\newcommand{\eers}{\end{exercise}}
\newcommand{\eexs}{\end{example}}
\newcommand{\eems}{\end{example}}
\newcommand{\els}{\end{lemma}}
\newcommand{\eles}{\end{lemmaex}}
\newcommand{\ets}{\end{theorem}}
\newcommand{\eds}{\end{definition}}
\newcommand{\eps}{\end{proposition}}

\newcommand{\bi}{\begin{itemize}}
\newcommand{\ei}{\end{itemize}}
\newcommand{\ben}{\begin{enumerate}}
\newcommand{\een}{\end{enumerate}}

\def\vbar{\mathchoice{\vrule height6.3ptdepth-.5ptwidth.8pt\kern-.8pt}
   {\vrule height6.3ptdepth-.5ptwidth.8pt\kern-.8pt}
   {\vrule height4.1ptdepth-.35ptwidth.6pt\kern-.6pt}
   {\vrule height3.1ptdepth-.25ptwidth.5pt\kern-.5pt}}
\def\fudge{\mathchoice{}{}{\mkern.5mu}{\mkern.8mu}}
\def\bbc#1#2{{\rm \mkern#2mu\vbar\mkern-#2mu#1}}
\def\bbb#1{{\rm I\mkern-3.5mu #1}}
\def\bba#1#2{{\rm #1\mkern-#2mu\fudge #1}}
\def\bb#1{{\count4=`#1 \advance\count4by-64 \ifcase\count4\or\bba A{11.5}\or
   \bbb B\or\bbc C{5}\or\bbb D\or\bbb E\or\bbb F \or\bbc G{5}\or\bbb H\or
   \bbb I\or\bbc J{3}\or\bbb K\or\bbb L \or\bbb M\or\bbb N\or\bbc O{5} \or
   \bbb P\or\bbc Q{5}\or\bbb R\or\bbc S{4.2}\or\bba T{10.5}\or\bbc U{5}\or
   \bba V{12}\or\bba W{16.5}\or\bba X{11}\or\bba Y{11.7}\or\bba Z{7.5}\fi}}

\def \qed {{\hspace*{\fill}$\halmos$\medskip}}

\def \A {{\cal{A}}}
\def \H {{\cal{H}}}

\def \G {{\cal{G}}}
\def \C {{\cal{C}}}
\def \D {{\cal{D}}}

\newcommand{\ba}[1]{\addtocounter{for}{1} \begin{eqnarray}\label{#1}}
\newcommand{\ea}{\end{eqnarray}}

\def\sqr#1#2{{\vcenter{\vbox{\hrule height .#2pt
                             \hbox{\vrule width .#2pt height#1pt \kern#1pt
                                   \vrule width .#2pt}
                             \hrule height .#2pt}}}}

\def\pmb#1{\setbox0=\hbox{#1}%
   \kern-.025em\copy0\kern-\wd0
   \kern.05em\copy0\kern-\wd0
 \kern-.025em\raise.0433em\box0 }
\def\sqr#1#2{{\vcenter{\vbox{\hrule height.#2pt
     \hbox{\vrule width.#2pt height#1pt \kern#1pt
   \vrule width.#2pt}\hrule height.#2pt}}}}

\def\B{{\mathcal B}}
\def\E{{\mathbb E}_0}
\def\P{{\mathbb P}}
\def \RR{{\mathcal R}}
\def \NN{{\mathcal N}}

\def\N{{\mathbb N}}
\def\Z{{\mathbb Z}}
\def\R{{\mathbb R}}

\def\var{\text{var}}

\def\reff#1{(\ref{#1})}

\def \ind {\hbox{ 1\hskip -3pt I}}
\newcommand {\cro}[1] {\left[ {#1} \right]}
\newcommand {\acc}[1] {\left\{ {#1} \right\}}
\newcommand {\pare}[1] {\left( {#1} \right)}

\newcommand {\bra}[1] {\left< {#1} \right>}
\newcommand {\sous}[1] {\underline{#1}}

\begin{document}

\title{Annealed lower tails for the energy of a polymer.}
\author{Amine Asselah \\ Universit\'e Paris-Est\\
Laboratoire d'Analyse et de Math\'ematiques Appliqu\'ees\\
UMR CNRS 8050\\ amine.asselah@univ-paris12.fr}
\date{}
\maketitle
\begin{abstract}
We consider the energy of a randomly charged polymer.
We assume that only charges on the same site interact pairwise.
We study the lower tails of the energy,
when  averaged over both randomness, in dimension three or more.
As a corollary, we obtain the {\it correct} temperature-scale
for the Gibbs measure.
\end{abstract}



{\em Running head}: Lower tails for the energy a polymer.

\section{Introduction} \label{sec-intro}
In this paper, we study the lower tails for the energy of a polymer.
This complements a companion paper \cite{A09upper}
dealing with the upper tails.
Lower and upper tails are different stories, and the two
papers are independent from each other, though they use the
same model, and the same notations. Thus, our polymer 
is a linear chain
of $n$ monomers each carrying a random charge, and sitting sequentially on
the positions of a symmetric random walk. 
\begin{itemize}
\item[(i)] The symmetric
random walk on $\Z^d$ is denoted $\{S(n),n\in \N\}$.
When $S(0)=z\in \Z^d$, its law is denoted ${\mathbb P}_z$.
\item[(ii)] The random field of charges is denoted $\{\eta(n), n\in \N\}$.
The charges are centered i.i.d.\,with a finite forth moment.
We denote by $\eta$ a generic charge variable, and the charges' law is
denoted by $Q$. 
\end{itemize}
The monomers interact pairwise only when they occupy the same site 
on the lattice. The interaction produces an energy
\be{intro.1}
H_n=\sum_{z\in \Z^d}
\sum_{0\le i\not= j< n} \eta(i)\eta(j) \ind\acc{S(i)=S(j)=z}.
\ee

Our toy-model comes from physics, where it is used to model
proteins or DNA {\it folding}. However, physicists' usual setting
differs from ours by three main features. (i) Their
polymer is usually {\it quenched}: a typical realization
of the charges is fixed, and the average is over the walk. 
(ii) A short-range repulsion is included by considering
random walks such as the self-avoiding walk or
the directed walk. (iii)
The averages are performed with respect to the
the so-called Gibbs measure: a probability measure obtained 
from $\P_0$ by weighting it with $\exp(\beta H_n)$,
with real parameter $\beta$. When $\beta$ is positive, the Gibbs
measure favors configuration with large energy; in other words, 
alike charges attract each other: this models {\it hydrophobic
interactions}, where the effect of avoiding the water solvent is mimicked
by an attraction among hydrophobic monomers. When $\beta$ is
negative, alike charges repel: this models Coulomb potential, and
describes also the effective repulsion between identical bases of RNA.
The issue is whether there is a {\it critical value} $\beta_c(n)$, such
that as $\beta$ crosses $\beta_c(n)$, a phase transition occurs.
For instance, Garel and Orland \cite{GO}
observed a phase transition as $\beta$ crosses a $\beta_c(n)\sim 1/n$, from
a collapsed shape to a random-walk like shape.
Kantor and Kardar~\cite{KK} discussed the quenched
model for the case $\beta<0$,
that is when alike charges repel. Some heuristics
(dimensional analysis on the continuum version) suggests that the
{\it (upper) critical} dimension is 2: for $d\ge 3$, the 
polymer looks like a simple random walk, whereas when $d<2$, its
average end-to-end distance is $n^\nu$ with $\nu=\frac{2}{d+2}$.
Let us also mention studies of
Derrida, Griffiths and Higgs \cite{dgh} and 
Derrida and Higgs \cite{dh}: both study the quenched
Gibbs measure $\exp(-\beta H_n)d\tilde \P_0$, with $\beta>0$,
for a one dimensional directed random walk $\tilde \P_0$, 
and obtain evidence for a phase transition (a so-called
weak freezing transition). 

Our interest stems from recent mathematical works of Chen~\cite{chen07}, and
Chen and Khoshnevisan~\cite{chen-khosh}, dealing with central limit
theorems for $H_n$. Chen~\cite{chen07} establishes also an annealed 
moderate deviation principle, under the additional assumption that 
$E[\exp(\lambda \eta^2)]<\infty$, for some $\lambda>0$. 
More precisely, with the annealed law denoted $P$, $d\ge 3$, 
$n^{\frac{1}{2}}\ll\sqrt{n}\xi_n\ll n^{\frac{2}{3}}$,
(for two positive sequences $\{a_n,b_n,n\in \N\}$,
we say that $a_n\ll b_n$, when $\limsup \frac{\log(a_n)}{\log(b_n)}<1$),
X.Chen has obtained
\be{intro.4}
\lim_{n\to \infty} \frac{1}{\xi_n^2} 
\log\pare{ P(\pm \frac{H_n}{\sqrt n}\ge \xi_n)}= -\frac{1}{2c_d},
\quad\text{where}\quad c_d=\sum_{n\ge 1} \P_0(S(n)=0).
\ee
Our study complements the work \cite{chen07}.
We study the annealed probability
that $\{- H_n> \xi_n\}$ for $\xi_n\ge n^{\frac{2}{3}}$.
Also, we consider the simplest aperiodic walk: the walk jumps to
a nearest neighbor site or stays still with equal probability.

As in \cite{A09upper}, 
we rewrite the energy into a convenient form. For $z\in \Z^d$,
and $n\in \N$, we call $l_n(z)$ the {\it local times},
and $\check{q}_n(z)$ the {\it local charges}. That is
\[
l_n(z)=\sum_{k=0}^{n-1} \ind\acc{S(k)=z},
\quad\text{and}\quad
\check{q}_n(z)=\sum_{k=0}^{n-1} \eta(k) \ind\acc{S(k)=z}.
\]
We write $H_n=\sum_{z}\check{X}_n(z)+Y_n(z)$ with
\[
\check{X}_n(z)=\check{q}_n^2(z)-l_n(z),
\quad \text{and}\quad 
Y_n(z)=l_n(z)-\sum_{i=0}^{n-1}\!\eta(k)^2\ind\acc{S(k)=z}.
\]
Now, 
\be{intro.10}
Y_n=\sum_{z\in\Z^d} Y_n(z)=\sum_{i=0}^{n-1}\pare{1-\eta^2(i)},
\ee
is a sum of centered independent random variables, and its
large deviation asymptotic are well known (see below
Remark~\ref{rem-Hlower}).
Thus, we focus on $\check{X}_n=\sum_{\Z^d} \check{X}_n(z)$. 

Before presenting our lower tails estimates, we provide
some heuristics.
\paragraph{Heuristics.} Since we are interested in annealed
estimates, note that
\be{heur-1}
\check{X}_n
\stackrel{\text{law}}{=} X_n:=\sum_{z\in \Z^d} 
l_n(z)\pare{\zeta_z(l_n(z))-1},\quad\text{where}\quad
\zeta_z(n)=\pare{\frac{1}{\sqrt n}\sum_{i=1}^n \eta_z(i)}^2,
\ee
where $\{\eta_z(i),z\in\Z^d,\ i\in \N\}$ is an i.i.d.\,sequence
with $\eta_z(i)\sim\eta$, and we still denote its law with $Q$.
Let us fix two lengths $T_n$ and $r_n$, and an energy $x_n$,
and estimate the cost
of folding $T_n$-monomers in a ball of radius $r_n$, say $B(r_n)$,
in order to realize 
\[
\sum_{z\in B(r_n)} l_n(z)\pare{1-\zeta_z(l_n(z))}\ge x_n.
\]
Note that necessarily $T_n\ge x_n$. Assume also that $T_n\gg |B(r_n)|$,
so that we expect many monomers to pile up on each site
of $B(r_n)$, and we further assume that the filling is uniform, that
is
\[
\forall z\in B(r_n),\quad l_n(z)\sim\frac{T_n}{|B(r_n)|}.
\]
Then, the optimal scenario comes up as we equate the cost
of the two constrains we are imposing. (i) We localize the walk
a time $T_n$ in a ball $B(r_n)$. This costs of the order of
$\exp(-\kappa T_n |B(r_n)|^{-2/d})$. (ii) We require the charges
to realize 
\be{heur-2}
\acc{\sum_{z\in B(r_n)} 1-\zeta_z(l_n(z))\ge \frac{x_n |B(r_n)|}{T_n}}.
\ee
Since, when we freeze the walk, the variables $\{1-\zeta_z(l_n(z)),
z\in B(r_n)\}$ are independent, centered and with finite variance
(if $E[\eta^4]<\infty$), the cost of \reff{heur-2} is
\be{heur-3}
P\pare{\sum_{z\in B(r_n)} 1-\zeta_z(l_n(z))\ge \frac{x_n |B(r_n)|}{T_n}}
\sim\exp\pare{ -\frac{x_n^2 |B(r_n)|}{T_n^2}}.
\ee
As we equate the two costs, we find
\be{heur-4}
\frac{x_n^2 |B(r_n)|}{T_n^2}=\frac{T_n}{|B(r_n)|^{2/d}}
\Longrightarrow |B(r_n)|^{\frac{d+2}{d}}=\frac{T_n^3}{x_n^2}.
\ee
Thus, the heuristic discussion suggests that for some constant $c>0$
\be{heur-5}
P(X_n\le -x_n)\sim 
\exp\pare{-c x_n^{\frac{4}{d+2}}T_n^{\frac{d-4}{d+2}}}.
\ee
Note that the exponent $\frac{d-4}{d+2}$ of $T_n$ in \reff{heur-5}
suggests that $d=3$ and $d>4$ have a distinct phenomenology. When 
$d=3$, the cheapest cost is reached when $T_n=n$: the polymer
is {\it entirely folded} in a ball of volume $(\frac{n^3}{x_n^2})
^{\frac{3}{5}}$. 
Also, the sum of local charges, $\check{q}_n$, 
over this domain performs a {\it moderate deviations}.

When $d>4$, the cheapest cost requires the
smallest $T_n$, which is $x_n\le n$. Thus, the polymer is
{\it partially folded}, and \reff{heur-5} implies that
the volume of the ball is $x_n^{\frac{d}{d+2}}$.
Also, on each site the {\it local charge} 
performs a {\it typical fluctuation}.

Our heuristics set the stage for the following mathematical statements.
\bt{prop-neg3}
Assume $d=3$, and $E[\eta^4]<\infty$.
There are constants $a_0,c_3^\pm$ such that for $a_0\le \xi_n<n^{1/3}$,
\be{neg.3}
\exp\pare{-c_3^- \xi_n^{\frac{4}{5}} n^{\frac{1}{3}}}
\le P(\check{X}_n\le -\xi_n n^{2/3})\le 
\exp\pare{-c_3^+ \xi_n^{\frac{4}{5}} n^{\frac{1}{3}}}.
\ee
Moreover, we have the following description of the dominant
strategy. For a constant $a$ large enough,
\be{strat-LT3}
\lim_{n\to\infty} P\pare{ 
|\{z\in \Z^d: \frac{ \xi_n^{\frac{6}{5}}}{a}\le l_n(z)\le
 a \xi_n^{\frac{6}{5}}\}|\ge  \frac{n}{a^4 \xi_n^{6/5}}\  \bigg\|\ 
\check{X}_n\le -\xi_n n^{\frac{2}{3}}}=1.
\ee
\et
In dimension 4 and more, there are two regimes. In the following
regime, the energy has the same behavior as in the 
moderate deviation regime, where the polymer is {\it unfolded}.
\bt{prop-neg4-I}
Assume $d\ge 4$, and $E[\eta^4]<\infty$.
For any $\epsilon$ positive, choose any sequence $\{\xi_n\}$ with
\[
\xi_n\in [n^{1/6},n^{(d/2)/(d+4)-\epsilon}].
\]
There are $c_1,c_2>0$, such that for $n$ large enough
\be{neg.1}
\exp\pare{ -c_1 \xi_n^2 }\le P\pare{ \check{X}_n\le -\xi_n\sqrt n}
\le \exp\pare{ -c_2 \xi_n^2}.
\ee
Moreover, for a constant $A$ large enough
\be{neg.1-bis}
\lim_{n\to\infty} P\pare{
\sum_{z:\  l_n(z)\ge A}\check{X}_n(z)\le -\xi_n \sqrt n}=0.
\ee
\et

The second regime corresponds to a {\it partially folded polymer} as
alluded to in the heuristic discussion.
\bt{prop-neg4-II}
Assume $d\ge 4$, and $n^{\frac{d+2}{d+4}}<\xi_n\le \xi n$ with $\xi<1$.
For a constant $c_d^-$, and for any $\epsilon>0$,
\be{neg.2}
\exp\pare{-c_d^-\xi_n^{\frac{d}{d+2}}}\le P\pare{ \check{X}_n\le -\xi_n}
\le \exp\pare{-\xi_n^{\frac{d}{d+2}}n^{-\epsilon}}.
\ee
\et
\br{rem-Hlower}
The lower tail behavior of $H_n$ depends on a competition
between $\check{X}_n$ and $Y_n$ whose upper tail behavior is given
in Remark~\ref{rem-Y}. Let us mention that if $\alpha\ge \frac{2d}{d+2}$,
then the lower tails of $H_n$ are identical to that of $\check{X}_n$.
When $d\ge 4$, and $\alpha< \frac{2d}{d+2}$, then $Y_n$ dictates
the behavior of $H_n$: the correct speed for the lower tails
of $H_n$ is $\min(\xi_n^2/n,\ \xi_n^{\alpha/2})$. In $d=3$, the
correct speed for the lower tails
of $H_n$ is $\min(\xi_n^{4/5}n^{-1/5},\ \xi_n^{\alpha/2})$.
Thus, as soon as $\alpha\ge 2$, 
the lower tails of $H_n$ are identical to that of $\check{X}_n$.
\er
\br{rem-weak}
The weakness in the upper bound in \reff{neg.2}
(the artifact $n^{-\epsilon}$ in the exponent) reflects a deep 
technical gap in estimating the distribution of the size of level
sets of the local times of the random walk. We state it as a conjecture.
\bcon{conj-1}
Assume $d\ge 3$, and let $\{y_n,n\in \N\}$ be
a sequence going to infinity, with $y_n^{1+d/2}\le n$.
Then, there is $\kappa_d>0$ (independent on $n$) such that
\be{eq-con1}
\P_0\pare{| \acc{z:\ l_n(z)\ge y_n}|\ge y_n^{d/2}}\le
\exp(-\kappa_d y_n^{d/2}).
\ee
\econ
One way to understand the difficulty of \reff{eq-con1} is
to see that the number of possible regions of volume
$y_n^{d/2}$ inside $[-n,n]^d$ exceeds $\exp(\kappa y_n^{d/2})$, for any
$\kappa>0$.
\er

We give now an elementary application of Theorem~\ref{prop-neg3}
to the study of annealed Gibbs measure in dimension three.
For simplicity, we further assume that $\eta\in \{-1,1\}$, so that
$H_n=\check{X}_n$. The annealed Gibbs measure is the following
probability measure: for $\beta>0$, we set
\be{def-gibbs}
dP_{n,\beta}^{-}= \frac{\exp(-\beta H_n)dP}{Z_n^-(\beta)}
\quad\text{where}\quad Z_n^-(\beta)=E\cro{\exp(-\beta H_n)}.
\ee
The normalizing constant $Z_n^-(\beta)$ is called partition function.
The measure $P_{n,\beta}^{-}$ favors configurations with large values of
$-H_n$, so that it forces local charges to neutralize.
When dealing with the Gibbs measure, the issue is to find the
{\it correct} temperature-scaling for which a phase-transition occurs.
Indeed, the interesting biological phenomenon which motivates
polymer modelling is {\it folding}, that is the process of going from
a (transient) random-walk shape to a globular-looking shape, under the
tuning of temperature, or salt-concentration. Thus, we expect
a critical parameter $\beta_c(n)$ (which might scale with the polymer
size), such that for $\beta>\beta_c(n)$, typical polymers are 
globular-like looking, whereas when $\beta<\beta_c(n)$, typical
polymers look like typical random walk trajectories.

Biskup and K\"onig \cite{biskup-konig}
(see also Buffet and Pul\'e \cite{buffet-pule})
obtain results and some heuristics 
on the {\it annealed} Gibbs measure (i.e.\,averaged over both randomness). 
They use that when freezing the random walk, and averaging over charges 
\be{biskup-konig1}
E_Q[e^{-\beta H_n}]=c_n
\exp(-\sum_{z\in \Z^d} V(l_n(z)))\quad\text{where for $x$ large}\quad 
V(x)\sim\frac{1}{2} \log(1+2\beta x),
\ee 
where $\beta>0$ and $c_n$ is a constant. When we assume 
that $Q(\eta=\pm 1)=\frac{1}{2}$, then $c_n=\exp(\beta n)$, and
the study \cite{biskup-konig}
suggests that when performing a further random walk average
\be{biskup-konig2}
e^{-\beta n}Z_n^-(\beta)=E\cro{e^{-\beta (H_n+n)}}
=\exp\pare{-\beta\chi n^{\frac{d}{d+2}}
\log(n)^{\frac{2}{d+2}}(1+o(1))}.
\ee 
and $\chi>0$ is independent of $\beta$.
Also, the proof of \cite{biskup-konig} suggests that, under the
annealed measure, the walk is localized a time $n$ into a ball
of volume $(n/\log(n))^{\frac{d}{d+2}}$.

Our results focus on determining the correct temperature-scale,
and are as follows.
\bp{prop-gibbs}
Assume that $d=3$, and $Q(\eta=\pm 1)=\frac{1}{2}$. 
The {\it correct} temperature-scaling is $1/n^{2/5}$.
More precisely, there are positive
constants $\beta_1<\beta_2$, and the following holds. When
$\beta>\beta_2$ (the {\it low temperature regime}), then for some
positive constants $a,c_1$
\be{gibbs-1}
\exp(\beta n^{3/5})\ge Z_n^-\pare{\frac{\beta}{n^{2/5}}}\ge
\exp(c_1\beta n^{3/5}),
\ee
and,
\be{gibbs-2}
\lim_{n\to\infty} P_{n,\frac{\beta}{n^{2/5}}}^{-}
\pare{|\{z\in \Z^d: \frac{ n^{\frac{2}{5}}}{a}\le l_n(z)\le
 a n^{\frac{2}{5}}\}|\ge  \frac{n^{3/5}}{a^4 }}=1.
\ee
When $\beta<\beta_1$ (the {\it high temperature regime}),
for $c_d$ defined in \reff{intro.4},
\be{gibbs-3}
\lim_{n\to\infty} \frac{1}{n^{1/5}}\log Z_n^-\pare{\frac{\beta}{n^{2/5}}}=
\ \frac{ c_d \beta^2}{2}.
\ee
Moreover, there is a positive constant $b$, such that
\be{gibbs-4}
\lim_{n\to\infty} P_{n,\frac{\beta}{n^{2/5}}}^{-}
\pare{\{z\in \Z^d:  l_n(z)\ge b n^{1/5} \}\not=\emptyset}=0.
\ee
\ep
\br{rem-gibbs} We stress that \reff{gibbs-4} is not the
`correct' result, since we expect that in the high temperature regime,
the polymer behaves like a random walk and we conjecture 
rather that for large $b$
\be{gibbs-5}
\lim_{n\to\infty} P_{n,\frac{\beta}{n^{2/5}}}^{-}\pare{\{z\in \Z^d:\
l_n(z)\ge b \log(n) \}\not=\emptyset}=0.
\ee
We include \reff{gibbs-4} to show the difference with \reff{gibbs-2}
which occurs in the low temperature regime.
\er

The paper is organized as follows. In Section~\ref{sec-prel}, we 
recall the large deviations for the $q$-norm of the local times. 
We have then divided Theorems~\ref{prop-neg3}, \ref{prop-neg4-I},
and \ref{prop-neg4-II}, into their upper bounds parts, and their
lower bounds parts. Upper bounds are treated in Section~\ref{sec-ublt},
while lower bounds are treated in Section~\ref{sec-lblt}.
Finally, Section~\ref{sec-gibbs} contains
the proof of Proposition~\ref{prop-gibbs}.

\section{Preliminaries}\label{sec-prel}
\subsection{Sums of Independent variables}
A. Nagaev has considered in \cite{anagaev} a sequence $\{\bar Y_n,n\in \N\}$ 
of independent centered i.i.d satisfying $\H_\alpha$ with $0<\alpha<1$, and
has obtained the following upper bound (see also inequality (2.32) of 
S.Nagaev~\cite{snagaev1}).
\bp{prop-anagaev} Assume $E[\bar Y_i]=0$ and
$E[(\bar Y_i)^2]\le 1$. There is a constant $C_Y$, such that 
for any integer $n$ and any positive $t$
\be{prel.4}
P\pare{\bar Y_1+\dots+\bar Y_n\ge t}\le C_Y\pare{
nP\pare{\bar Y_1>\frac{t}{2}}+ \exp\pare{-\frac{t^2}{20 n}} }.
\ee
\ep
\br{rem-Y} 
Note that if $\eta\in \H_\alpha$ for $1<\alpha\le 2$, 
then $\eta^2\in \H_{\frac{\alpha}{2}}$. Thus, for $\bar Y_i=\eta(i)^2-1$,
Proposition~\ref{prop-anagaev} yields
\be{term-Y1}
P\pare{\sum_{i=1}^n (\eta(i)^2-1)\ge \xi_n}\le 
C_Y\pare{n\exp\pare{-c_{\alpha}(\xi_n)^{\alpha/2}}+
\exp\pare{ -\frac{\xi^2 n^{2\beta-1}}{20}}}.
\ee
\er

Finally, we specialize to our setting a general lower bound 
of S.Nagaev (see Theorem 1 of~\cite{snagaev2}). 
Let $\{\Lambda_n,\ n\in \N\}$
a sequence of subsets of $\Z^d$, and for each
$n$, let $\{Y_z^{(n)},z\in \Lambda_n\}$ be independent and
centered random variables. Let 
\[
\sigma_n^2=\sum_{z\in \Lambda_n} E\cro{(Y_z^{(n)})^2},\quad
\text{and}\quad \C_n^3=\sum_{z\in \Lambda_n} E\cro{|Y_z^{(n)}|^3}.
\]

\bp{prop-snagaev} Consider a sequence $\{t_n,\ n\in \N\}$ such that
for a small enough $\epsilon_\NN>0$ and $n$ large enough
\be{cond-nagaev}
1\le t_n\le\ \epsilon_\NN\ \min( \frac{\sigma_n^3}{\C_n^3},
\sigma_n (\max_{z\in \Lambda_n} \sqrt{ E[(Y_z^{(n)})^2]})^{-1}\ ),
\ee
then, there is a positive constant $\kappa$ such that
\be{gauss.13}
P\pare{ \frac{1}{\sigma_n}\sum_{z\in \Lambda_n}Y_z^{(n)}\ge t_n}\ge
\exp\pare{-\frac{t_n^2}{2}(1+\epsilon_\NN\kappa)}.
\ee
\ep
\subsection{On self-intersection local times}\label{sec-silt}
In this section, we recall and establish useful estimates
for functionals of the local times. First, for
any $z\in \Z^d$, we estimate the variance of $q_n^2(z)-l_n(z)$
\be{LB.9}
q_n^2(z)-l_n(z) =\big(\sum_{i\le l_n(z)}\eta_z(i)\big)^2-l_n(z)
=\sum_{i\le l_n(z)}(\eta^2_z(i)-1)
+2\sum_{1\le i<j\le l_n(z)} \eta_z(i)\eta_z(j),
\ee
It is immediate to obtain, for $\chi_1=E[\eta^4]+1$
\be{LB.10}
2 \pare{l_n^2(z)-l_n(z)}\le E_Q\cro{(q_n^2(z)-l_n(z))^2}
= l_n(z)\pare{E_Q[\eta^4]-1}+2\pare{l_n^2(z)-l_n(z)}\le
\chi_1 l_n^2(z).
\ee
Second, we summarize the asymptotic behavior of
the $q$-norm of local times (for any real $q>1$)
\be{def-lq}
\|l_n\|_q^q=\sum_{z\in \Z^d} l_n^q(z).
\ee 
In dimension three and more, Becker and K\"onig \cite{BK} have shown that
there are positive constants, say $\kappa(q,d)$, such that almost surely
\be{BK-main}
\lim_{n\to\infty} \frac{\|l_n\|_q^q}{n}=\kappa(q,d).
\ee
The large deviations, and central limit theorem for $\|l_n\|_q$
are tackled in \cite{A08}: we establish a shape transition 
in the walk's strategy to realize
the deviations $\{\|l_n\|_q^q-E[\|l_n\|_q^q]\ge n\xi\}$ with $\xi>0$.
This transition occurs at a critical 
value $q_c(d)=\frac{d}{d-2}$ suggesting the following picture.
\begin{itemize}
\item In the {\it super-critical regime} $q>q_c(d)$, 
the walk performs a short-time clumping on finitely many sites.
\item In the {\it sub-critical regime} $q<q_c(d)$, 
the walk is localized during the whole
time-period in a ball of volume $n/\xi^{\frac{1}{q-1}}$ where
it visits each site of the order of $\xi^{\frac{1}{q-1}}$-times.
\end{itemize}
We first recall Theorem 1.2 of \cite{A08} which deals with the 
{\it super-critical} regime. 
\bl{theo-sup} Assume $d\ge 3$ and $q>q_c(d)$. There are constants 
$C,c(q,d)$
(depending only on $d$ and $q$), such that for $\xi_n\ge 1$, and any
integer $n$ 
\be{eq-sup1}
\P_0\pare{ \|l_n\|_q^q-\E\cro{\|l_n\|_q^q}>\xi_n n}\le 
C \exp\pare{ -c(q,d)(\xi_n\ n)^{\frac{1}{q}}}.
\ee
\el
Also, Lemma 1.4 of \cite{A08} estimates the cost of the contribution of
{\it low } level sets to an excess $q$-norm. Thus, define for $x,y>0$
\[
\D_n(x,y):=\acc{z:\ x<l_n(z)\le y}.
\]
\bl{lem-sup}
Assume $d\ge 3$ and $q\ge q_c(d)$. For $\gamma\ge 1$, and $\chi>0$
and $\epsilon>0$, there is a constant $C$ such that for any sequence $y_n$ 
\be{super.1}
\P_0\pare{ \sum_{z\in \sous\D_n(1,y_n)} l_n^q(z)\ge \chi n^\gamma}\le C
\exp\pare{ -\frac{n^{\gamma/q_c(d)-\epsilon}}{y_n^{(q/q_c(d)-1)}}}.
\ee
When $\gamma=1$, one needs to take $\chi>\kappa(q,d)$ in \reff{super.1}.
\el
\br{rem-delicate}
Actually Lemma 1.4 of \cite{A08} is only stated for $\gamma>1$. 
An inspection of its proof, shows that it covers also the case
$\gamma=1$ provided that $\chi>\kappa(q,d)$.
In \reff{super.1}, we are unable to get rid of the $\epsilon$. 
This is a delicate
issue which is also responsible for a gap in the exponent of 
the speed in Region III of \cite{AC05} (inequality (8)).
\er
The next result deals with {\it sub-critical regime}.
It follows from Theorem 1.1 and Remark 1.3 of \cite{A08}.
\bl{lem-sub} Assume $d\ge 3$ and $1<q<q_c(d)$. There are constants 
$C,c(q,d)$
(depending only on $d$ and $q$), such that for $\xi_n\ge 1$, and any
integer $n$ 
\be{append.1}
\P_0\pare{ \|l_n\|_q^q-\E\cro{\|l_n\|_q^q}>\xi_n n}\le 
C \exp\pare{ -c(q,d)\xi_n^{\frac{2}{d}\frac{1}{q-1}}
n^{1-\frac{2}{d}}}.
\ee
\el
\br{rem-mistake}
For $d=3$, \reff{append.1} is mistakenly reported in
\cite{A06}. Fortunately, this is of no consequence since
(with the notations of \cite{A06} and in the so-called Region II), we need
there
\[
\frac{2}{3}(\beta+b)-\frac{1}{3}-\epsilon> \beta-b\Longleftrightarrow
5\frac{\beta}{\alpha+1}> \beta+1+3\epsilon\Longleftrightarrow
\beta>\frac{\alpha+1}{4-\alpha}.
\]
This latter condition defines Region II.
\er
We now state a corollary of Lemmas~\ref{lem-sup} and \ref{lem-sub}, 
whose immediate proof is omitted.
\bc{cor-2beta}
Assume $d\ge 3$ and $\xi_n\ge n^{\frac{2}{3}}$. 
For $\epsilon>0$ small enough, and $n$ large enough
\be{neglect-2beta}
\P_0\pare{ \|l_n\|_2 \ge \xi_n n^{-\epsilon}}\le
\exp\pare{ -\xi_n^{\frac{d}{d+2}}n^{\epsilon}}.
\ee
\ec
\section{Upper Bounds.}
\label{sec-ublt}
In this section, we prove the upper bounds in
Theorems \ref{prop-neg3}, \ref{prop-neg4-I}, and \ref{prop-neg4-II}.
When dealing with large deviations, a natural approach 
is to perform a Chebychev's exponential inequality. 
If we expect $P(X_n\le -x_n)\sim\exp(-\zeta_n)$,
then for $\lambda>0$, and $y_n=x_n/\zeta_n$
\be{neg.5}
P\pare{\bra{l_n,1-\zeta_.(l_n)}\ge x_n}\le
e^{-\lambda \zeta_n}\ E\cro{ \exp\pare{\lambda
\bra{\frac{l_n}{y_n},1-\zeta_.(l_n)}}}.
\ee
Now, to get rid of the dependence between field and local time,
we first perform an integration over the charges. We define
for $x\in \R^+$ and $n\in \N$
\be{neg.6}
\tilde\Gamma(x,n)=\log E_Q\cro{\exp\pare{x (1-\zeta_0(n))}}.
\ee
Since $1-\zeta_0(n)\le 1$, and since $e^u\le 1+u+u^2$ when $u\le 1$,
we have, for the constant $\chi_1$ which appears in \reff{LB.10},
\be{neg.7}
\begin{split}
\tilde\Gamma(x,n)
\le &\ind_{\{x\ge 1\}} x+\ind_{\{x< 1\}}
\ \log E_Q\cro{ 1+ x(1-\zeta_0(n))+ x^2 (1-\zeta_0(n))^2}\\
\le & \ind_{\{x\ge 1\}}x+\ind_{\{x< 1\}} \log\pare{1+x^2 \var(\zeta_0(n))}\\
\le & \ind_{\{x\ge 1\}}x+\ind_{\{x< 1\}} x^2 \sup_k \var(\zeta_0(k))
\le \ind_{\{x\ge 1\}}x+\ind_{\{x< 1\}}\ \chi_1 x^2.
\end{split}
\ee
\br{rem-hoffding}
Note first that \reff{neg.7} implies that $\tilde\Gamma(x,n)\le
\max(1,\chi_1) x^2$. Secondly, the dependence of $\tilde\Gamma(x,n)$ 
on the local times has vanished in these two regimes.
\er
Using \reff{neg.5} and \reff{neg.6}, our first step is
\be{neg.main}
P\pare{\bra{l_n,1-\zeta_.(l_n)}\ge x_n}\le
e^{-\lambda \zeta_n}\ \E\cro{\exp\pare{
\sum_{z\in \Z^d} \tilde\Gamma(\frac{\lambda l_n(z)}{y_n},l_n(z))}}.
\ee
We introduce some notations. For $0<x<y$, and $\chi>0$
\be{neg.8}
\D_n(x,y)=\acc{z\in \Z^d:\ x<l_n(z)\le y },\quad
\text{and}\quad 
\B(x,y;\chi)= \acc{\sum_{z\in\D_n(x,y)} l_n^2(z)\ge \chi}.
\ee
Also, we add a handy notations: for a subset $\Lambda\subset \Z^d$,
 $X_n(\Lambda)=\sum_{z\in \Lambda} X_n(z)$.

To treat separately the contribution of the two regimes of $\tilde\Gamma$,
we divide the visited sites of the walk into $\D_n(1,y_n)$,
and $\D_n(y_n,n)$. 
For $x_n'=x_n''=x_n/2$, and $0<\lambda<1$, we 
abbreviate $\B(1,y_n;\chi y_nx_n)$ by $\B$, and we have
\be{neg.10}
\begin{split}
P\big(-X_n\ge x_n\big)\le&
\P_0\pare{ l_n(\D_n(y_n,n))\ge x_n'}+
P\big(-X_n(\D_n(1,y_n))\ge x_n''\big)\\
\le &
\P_0\pare{ l_n(\D_n(y_n,n))\ge x_n'}+\P_0\pare{\B}+
P\big(-X_n(\D_n(1,y_n))\ge x_n'' ,\ \B^c\big)\\
\le &
\P_0\pare{ l_n(\D_n(y_n,n))\ge x_n'}+\P_0\pare{\B}\\
&\qquad\qquad
+\exp\pare{-\lambda \frac{x_n''}{y_n}}
\ \E\cro{\ind_{\B^c} \exp\pare{ \chi_1 \lambda^2
\sum_{\D_n(1,y_n)} \pare{\frac{l_n(z)}{y_n}}^2}}\\
\le & 
\P_0\pare{ l_n(\D_n(y_n,n))\ge x_n'}+\P_0\pare{\B}+
\exp\pare{-\zeta_n( \frac{\lambda}{2}-\lambda^2\chi_1\chi)}.
\end{split}
\ee
Note that the occurrence of an $l_2$-norm of the local time, in $\B(1,y_n;\chi)$,
is not arbitrary but is a consequence of the asymptotic of the
log-Laplace in \reff{neg.7}.

We discuss now the respective contributions of the {\it top level term}
$\{l_n(\D_n(y_n,n))\ge x_n'\}$, and 
of the {\it bottom level term} $\B(1,y_n;\chi y_nx_n)$.
Note that the threshold $y_n$ defining the {\it top level term}
is determined by the log-Laplace, and may not be the value
of the level set having a dominant contribution to our large deviation.
\paragraph{Top level term.}
First, note that for any $q>1$,
\be{neg.9}
\acc{l_n(\D_n(y_n,n))\ge x_n'}\subset
\acc{\|\ind_{\D_n(y_n,n)}l_n\|_q^q\ge \frac{1}{2} x_ny_n^{q-1}}.
\ee
The event on the right hand side of \reff{neg.9} has
a small probability if $x_ny_n^{q-1}>\kappa(q,d) n$, 
where $\kappa(q,d)$ is defined in \reff{BK-main}. 

We distinguish $q<q_c(d)$ and $q>q_c(d)$ with $q_c(d)=d/(d-2)$ (see
Section~\ref{sec-silt}). (i) When
$q<q_c(d)$, the so-called {\it subcritical regime},
Lemma~\ref{lem-sub} yields
\be{neg.30}
P\pare{\|\ind_{\D_n(y_n,n)}l_n\|_q^q\ge \frac{1}{2}x_ny_n^{q-1}}\le
\exp\pare{-c(q,d) \pare{ \frac{x_n}{2n} y_n^{q-1}}^{\frac{2}{d}
\frac{1}{(q-1)}} n^{1/q_c(d)}}.
\ee
Now, since $x_n\le n$, the map $q\mapsto \frac{x_n}{n}^{\frac{1}{(q-1)}}$
increases on $[1,q_c(d)[$. (ii) When $q>q_c(d)$, it is easy to check
that the upper bound given by Lemma~\ref{lem-sup},
increases on $]q_c(d),\infty[$, as a function of $q$. 
Thus, the best estimates we can
obtain on $\{l_n(\D_n(y_n,n))\ge x_n'\}$ is with a bound as in
\reff{neg.9} right at $q_c(d)$, for which we do not have sharp
estimates.

\paragraph{Bottom level term.}
When $2<q_c(d)$ (that is in $d=3$), we expect $\B(1,y_n;\chi y_nx_n)$
to be of order $\{\|l_n\|_2^2\ge \chi y_nx_n\}$, and by
Lemma~\ref{lem-sub}, we have in $d=3$, for $\chi x_n y_n>
\kappa(2,d) n$, that
\be{neg.31}
P(\B(1,y_n;\chi y_nx_n))\le 
P\pare{\|l_n\|_2^2\ge  \chi y_nx_n}\le \exp\pare{-c(2,3) (\chi y_nx_n)^{2/3}
n^{-1/3}}.
\ee
In this case, the cost of the bottom level set dominates the
top level sets, and it is therefore useless to consider $q>2$ in \reff{neg.30},
when $d=3$. When $q_c(d)\le 2$ (that is when $d\ge 4$), and $x_ny_n/n\to
\infty$, we can use Lemma~\ref{lem-sup}, even though this is not an optimal
result. 

It is clear from this discussion that the behavior of the lower tail
is distinct in $d=3$ and in $d\ge 4$. This leads to different
strategies, and different exponents. We discuss separately the case
$d=3$ and the case $d\ge 4$.
\subsection{Dimension 3}
We first make explicit the notations of \reff{neg.5}
\be{neg.32}
x_n=\xi_n n^{\frac{2}{3}},\quad \zeta_n=\xi_n^{\frac{4}{5}} n^{\frac{1}{3}},
\quad\text{and}\quad y_n=\frac{x_n}{\zeta_n}=\xi_n^{\frac{1}{5}} n^{\frac{1}{3}}.
\ee
where $\xi_n$ can vary in $[a_0,n^{\frac{1}{3}}]$, for a constant $a_0$ to
be specified later. Our first result is the following rough upper bound.
\bl{lem-roughUB} Assume $d=3$. There are positive constants $a_0,c_3^+$, such
that for $\xi_n\in [a_0,n^{\frac{1}{3}}]$ 
\be{neg.33}
P(-X_n\ge \xi_n n^{2/3})\le 3 
\exp\pare{-c_3^+ \xi_n^{\frac{4}{5}} n^{\frac{1}{3}}}.
\ee
\el
Note that in Section~\ref{sec-LBd3}, we establish a similar lower bound.

\noindent{\bf Proof of Lemma~\ref{lem-roughUB}}
Recall that \reff{neg.9}, for $q=2$, requires that $x_n y_n>2\kappa(2,3)n$,
which is equivalent to $\xi_n>a_0:=(2\kappa(2,3))^{5/6}$. Recall
that \reff{neg.31} requires that $\chi x_n y_n>\kappa(2,3) n$, which
is equivalent to $\chi \xi_n^{6/5}>\kappa(2,3)$, which in turn
requires that $\chi>1/2$.
Combining inequalities \reff{neg.10}, \reff{neg.9} with $q=2$, and 
\reff{neg.31}, we obtain for $0\le \lambda\le 1$
\be{neg.36}
P(-X_n\ge \xi_n n^{2/3})\le \exp\pare{ -\frac{c(2,3)}{2^{2/3}} \zeta_n}+
\exp\pare{ -c(2,3)\chi^{2/3} \zeta_n}+
\exp\pare{ -(\frac{\lambda}{2}-\lambda^2 \chi_1\chi) \zeta_n}.
\ee
We choose $\chi=1/4$, and $\lambda=\min(1/\chi_1,1)$ in \reff{neg.31}
to obtain the desired result.

\qed

\subsubsection{Upper bound in Theorem~\ref{prop-neg3}: $x_n=\xi_n n^{2/3}< n$}
We show in this section that the dominant {\it level set} of the local
times is of order $\xi_n^{\frac{6}{5}}$ much smaller than
$y_n$ when $x_n$ is much smaller than $n$.
We actually consider $x_n< a_1 n$ with $a_1$ to be chosen later small. 
For a large constant $a>0$, to be chosen later, 
we decompose $\{z:l_n(z)>0\}$ into $\D_1\cup\dots\cup\D_4$ with
\be{neg.34}
\D_1=\D_n(1,\frac{1}{a} \xi_n^{\frac{6}{5}}),\ 
\D_2=\D_n(\frac{1}{a} \xi_n^{\frac{6}{5}}, a \xi_n^{\frac{6}{5}}),\ 
\D_3=\D_n( a \xi_n^{\frac{6}{5}},\frac{y_n}{a}),\ \text{and}\quad
\D_4=\D_n(\frac{y_n}{a},n).
\ee
We then write
\be{neg.35}
P(-X_n\ge \xi_n n^{\frac{2}{3}})\le \sum_{i\not= 2}
P\pare{-X_n(\D_i)\ge \frac{1}{4} \xi_n n^{\frac{2}{3}}}+
P\pare{-X_n\ge \xi_n n^{\frac{2}{3}},\ 
-X_n(\D_2)\ge \frac{1}{4} \xi_n n^{\frac{2}{3}}}.
\ee
We now show that the contribution of $\D_2$ is the dominant one.

\noindent{\it a) Contribution of $\D_1$.}

We use Chebychev's inequality with $\lambda>0$,
\be{neg.37}
P\pare{-X_n(\D_1)\ge \frac{1}{4}\xi_n n^{2/3}}\le
e^{-\frac{\lambda}{4} \zeta_n} \E\cro{
\prod_{z\in \D_1} \exp\pare{\tilde \Gamma(\frac{\lambda l_n(z)}{y_n},l_n(z))}}.
\ee
Now, to justify the expansion of $\tilde \Gamma$ at 0, we need 
$\lambda \xi_n^{6/5}\le a y_n$ which is equivalent to $\lambda \xi_n\le 
an^{1/3}$. Assume that this latter fact holds. We have by \reff{neg.7}
\be{neg.38}
P\pare{-X_n(\D_1)\ge \frac{1}{4}\xi_n n^{2/3}}\le
\exp\pare{ -\frac{\lambda}{4}\zeta_n+\chi_1 \lambda^2 
\sum_{z\in \D_1} \frac{l_n^2(z)}{y_n^2}}.
\ee
It will be convenient to define $\chi_2=\max(\chi_1,\frac{1}{8})$.
We now use that $l_n(\D_1)\le n$, so that
\be{neg.39}
\sum_{z\in \D_1} \frac{l_n^2(z)}{y_n^2}\le 
\frac{\xi_n^{6/5}}{a y_n^2} l_n(\D_1) \le \frac{\xi_n^{6/5}n}{a y_n^2}=
\frac{\zeta_n}{a}.
\ee
We choose $\lambda=a/(8\chi_2)\le a n^{1/3}/\xi_n$ , and use \reff{neg.39}
in \reff{neg.38}
\be{neg.40}
P\pare{-X_n(\D_1)\ge \frac{1}{4}\xi_n n^{2/3}}\le
\exp\pare{ -\frac{a}{8^2\chi_2}\zeta_n}.
\ee

\noindent{\it b) Contribution of $\D_3$.}

For $0\le \lambda\le a$, and $\chi$ to be chosen later, we have
\be{neg.41}
\begin{split}
P\pare{-X_n(\D_3)\ge \frac{1}{4}\xi_n n^{2/3}}\le&
P\pare{ \B(a\xi_n^{6/5},y_n; \chi x_n y_n)}+
e^{-\frac{\lambda}{4} \zeta_n} \E\cro{\ind_{\B(.)^c}
\exp\pare{\chi_1 \lambda^2 \sum_{z\in \D_3} \frac{l_n^2(z)}{y_n^2}}}\\
\le & P\pare{ \B(a\xi_n^{6/5},y_n; \chi x_n y_n)}+
\exp\pare{ -(\frac{\lambda}{4}-\chi_1 \lambda^2\chi)\zeta_n}.
\end{split}
\ee
Choose $2<q<q_c(3)=3$, and by Lemma~\ref{lem-sub}
\be{neg.42}
\begin{split}
P\pare{ \B(a\xi_n^{6/5},y_n; \chi x_n y_n)}\le&
P\pare{\|l_n\|_q^q\ge (a\xi_n^{6/5})^{q-2} \chi x_n y_n}=
P\pare{\|l_n\|_q^q\ge a^{q-2}\xi_n^{6/5(q-1)} \chi n}\\
\le & \exp\pare{-c(q,3)\pare{ a^{q-2} \chi \xi_n^{\frac{6}{5}(q-1)}}^{
\frac{2}{3(q-1)}} n^{1/3}}\\
\le &  \exp\pare{-c(q,3)\pare{ a^{q-2} \chi}^{\frac{2}{3(q-1)}} \zeta_n}
\end{split}
\ee
Now, collecting \reff{neg.41} and \reff{neg.42}, we choose 
$\chi=a^{1-q/2}$ and for $a^{4-q}>(8\chi_1)^{-2}$ we have that
the optimal $\lambda$ in \reff{neg.41} satisfies $\lambda\le a$,
and
\be{neg.43}
\begin{split}
P\pare{-X_n(\D_3)\ge \frac{1}{4}\xi_n n^{2/3}}\le&
\exp\pare{-c(q,3)\pare{ a^{q-2} \chi}^{\frac{2}{3(q-1)}} \zeta_n}+
\exp\pare{ -(\frac{\lambda}{4}-\chi_1 \lambda^2\chi)\zeta_n}\\
\le & 
\exp\pare{-c(q,3)a^{\frac{q-2}{3(q-1)}} \zeta_n}+
\exp\pare{ -\frac{1}{8^2\chi_1}a^{q/2-1}\zeta_n}\\
\end{split}
\ee

\noindent{\it c) Contribution of $\D_4$.}

We proceed as in \reff{neg.9} and \reff{neg.30}.
\be{neg.44}
\begin{split}
P\pare{-X_n(\D_4)\ge \frac{1}{4}\xi_n n^{2/3}}\le&
P\pare{l_n(\D_4)\ge \frac{1}{4}\xi_n n^{2/3}}\le
P\pare{\|l_n\|_q^q\ge \frac{1}{4} \xi_n(\frac{y_n}{a})^{q-1}n^{2/3}}\\
\le &
\exp\pare{-c(q,3)\pare{\frac{\xi_n}{4n^{1/3}} 
(\frac{y_n}{a})^{q-1}}^{\frac{2}{3(q-1)}} n^{1/3}}.
\end{split}
\ee
Now, for $A>0$, and $2<q<3$,
\be{neg.45}
\frac{1}{a^{2/3}}(\xi_n y_n^{q-1})^{\frac{2}{3(q-1)}}\ 
 n^{\frac{1}{3}(1-\frac{2}{3(q-1))}}\ge A
\xi_n^{4/5} n^{1/3}\quad\Longleftrightarrow\quad
 \xi_n (aA^{3/2})^{\frac{(q-1)}{(q-2)}}\le n^{1/3}.
\ee
Our assumption is that $\xi_n<a_1 n^{1/3}$, and this implies that
\be{neg.46}
P\pare{-X_n(\D_4)\ge \frac{1}{4}\xi_n n^{2/3}}\le 
\exp\pare{-c(q,3) \frac{\zeta_n}{a_1^\gamma a^{2/3}}},\quad\text{with}\quad
\gamma=\frac{2(q-2)}{3(q-1)}>0.
\ee

\noindent{\it d) Contribution of $\D_2$.}

We recall the rough lower bound
$P(-X_n\ge \xi_n n^{\frac{2}{3}})\ge \exp\pare{- c_3^- \zeta_n}$,
and express \reff{neg.35} as
\be{neg.47}
P(-X_n\ge \xi_n n^{\frac{2}{3}})\le \sum_{i\not= 2}
P\pare{-X_n(\D_i)\ge \frac{1}{4} \xi_n n^{\frac{2}{3}}}+
P\pare{-X_n\ge \xi_n n^{\frac{2}{3}},
-X_n(\D_2)\ge\frac{1}{4} \xi_n n^{\frac{2}{3}}}.
\ee
When $a$ is large enough in \reff{neg.40} and \reff{neg.43}, and $a_1$
small enough in \reff{neg.46}, the terms with $\D_1$ and $\D_3$
are negligible. We then write
\be{neg.54}
\acc{-X_n(\D_2)\ge\frac{1}{4} \xi_n n^{\frac{2}{3}}}\subset
\acc{|\D_2|\ge \frac{n}{a^4 \xi_n^{6/5} }}\cup
\acc{\sum_{\D_2}\pare{1-\zeta_z(l_n(z))}\ge \frac{n^{\frac{2}{3}}}
{4 a \xi_n^{1/5}},\ |\D_2|\le \frac{n}{a^4 \xi_n^{6/5}}}.
\ee
Now, for dealing with the last event in \reff{neg.54}, note
that
\be{neg.55}
\acc{\sum_{\D_2}\pare{1-\zeta_z(l_n(z))}\ge \frac{n^{\frac{2}{3}}}
{4 a \xi_n^{1/5}},\ |\D_2|\le \frac{n}{a^4 \xi_n^{6/5}}}\subset
\acc{\frac{1}{\sqrt{|\D_2|}}\sum_{\D_2} \pare{1-\zeta_z(l_n(z))}\ge 
\frac{a \xi_n^{2/5} n^{1/6}}{4}}.
\ee
Now, we fix the randomness of the walk, and use that
$1-\zeta_z\le 1$, $E_Q[1-\zeta_z]=0$ and $E_Q[(1-\zeta_z)^2]\le \chi_1$
to obtain that (recall that $\zeta_n=\xi_n^{4/5} n^{1/3}$)
\be{neg.56}
P\pare{\frac{1}{\sqrt{|\D_2|}}\sum_{\D_2} \pare{1-\zeta_z(l_n(z))}\ge 
\frac{a \xi_n^{2/5} n^{1/6}}{4}}\le \exp(-\frac{a^2 \zeta_n}{4}).
\ee
We put together \reff{neg.47}, \reff{neg.54} and \reff{neg.56} to obtain
for $a$ large enough
\be{neg.49}
\lim_{n\to\infty} P\pare{ |\D_2|\ge \frac{n}{a^4 \xi_n^{6/5}} \ \bigg\|
\ -X_n\ge \xi_n n^{\frac{2}{3}}}=1.
\ee
\subsubsection{Upper bound in Theorem~\ref{prop-neg3}: $x_n=\xi n$ 
with $1>\xi>a_1$.}\label{sec-forgibbs}
Note that
\[
\xi_n=\xi n^{1/3},\quad
\zeta_n=\xi^{4/5} n^{3/5},\quad\text{and}\quad y_n=\xi^{1/5} n^{2/5}.
\]
Note that $\xi_n^{6/5}=\xi y_n<y_n$.
For a large constant $b>0$, to be
chosen later, we decompose $\{z:l_n(z)>0\}$ into $\D_1\cup\dots\cup\D_3$ with
\be{neg.50}
\D_1=\D_n(1,\frac{1}{b} \xi^{\frac{6}{5}}n^{2/5}),\
\D_2=\D_n(\frac{1}{b} \xi^{\frac{6}{5}}n^{2/5}, b\xi^{\frac{1}{5}}n^{2/5})
,\ \text{and}\quad\D_3=\D_n( by_n,n).
\ee
We then write
\be{neg.51}
P(-X_n\ge \xi n)\le \sum_{i\not= 2}P\pare{-X_n(\D_i)\ge \frac{1}{4} \xi n}+
P\pare{-X_n(\D_2)\ge \frac{1}{2} \xi n,\ -X_n\ge \xi n},
\ee
and we show that the contribution of $\D_2$ is the dominant one.

The treatment of $\D_1$ is similar to the previous case a). The choice
$\lambda=b/(8\chi_2)$ requires $\xi\le 8 \chi_2$, which holds since
$\xi<1\le 8 \chi_2$.

Then, for $\D_3$, we write
\be{neg.52}
\begin{split}
P\pare{-X_n(\D_3)\ge \frac{1}{4}\xi n}\le&
P\pare{l_n(\D_3)\ge \frac{1}{4}\xi n}\le
P\pare{\|l_n\|_2^2\ge \frac{1}{4} b\xi^{6/5} n^{2/5}n}\\
\le &
\exp\pare{-c(2,3)\pare{\frac{b}{4}}^{2/3}\zeta_n}.
\end{split}
\ee
By taking $b$ large enough, and proceeding as in the previous case d),
we reach that for $\xi<1$
\be{neg.53}
\lim_{n\to\infty} P\pare{ |\D_2|\ge \frac{\xi^{4/5}n^{3/5}}{b}\  \bigg\|\ 
\ X_n\le- \xi n }=1.
\ee
\subsection{Dimension 4 or more.}
We choose here $x_n,y_n$ and $\zeta_n$ as follows.
\be{choice-d4}
x_n=\xi_n \sqrt n,\quad \zeta_n=\xi_n^2,\quad\text{and}\quad
y_n=\frac{\sqrt n}{\xi_n}.
\ee
We first deal with the case $a_0 n^{1/6}\le
\xi_n\ll n^{\gamma_d -\epsilon}$, with $\gamma_d=(d/2)/(d+4)$,
and any $\epsilon$ positive.
\subsubsection{Proof of the Upper bound in \reff{neg.1}.}
Our starting point is the inequality \reff{neg.10} with
$x_n,y_n,\zeta_n$ as in \reff{choice-d4}. We deal
with each term on the right hand side of \reff{neg.10}.

First, choose $\chi>\kappa(2,d)$, and Lemma~\ref{lem-sup} gives
\be{neg.17}
P(\B(1,y_n;\chi x_n y_n))=\P_0\pare{
\sum_{z\in \D_n(1,y_n)} l_n^2(z)\ge \chi n}\le
\exp\pare{-\frac{n^{1/q_c(d)-\epsilon}}{y_n^{(2/q_c(d)-1)}}}.
\ee
Second, $n^{1/q_c(d)-\epsilon}\ge y_n^{(2/q_c(d)-1)} \xi_n^2$ is
equivalent to asking $\xi_n^{1+4/d}\le n^{1/2-\epsilon}$, which
is exactly the condition which defines this regime.

Now, we deal with the event $\{ l_n(\D_n(y_n, n))\ge x_n/2\}$. 
The proof of Proposition 3.3 of \cite{AC05} yields
\be{neg-prop3.3}
P\pare{  l_n(\D_n(y_n, n))\ge x_n/2}\le \exp\pare{-x_n^{1/q_c(d)}
y_n^{2/d}},
\ee
provided that for some fixed $a$ and $n$ large
\be{cond-top}
y_n^{1+\frac{2}{d}}\ge \log^a(n) x_n^{2/d}.
\ee
Now both $x_n^{1/q_c(d)}y_n^{2/d}\gg \xi_n^2$ and
condition \reff{cond-top} follow from
$\log\xi_n\le (d/2-\epsilon)/(d+4)\log(n)$.
Thus, for any $\epsilon>0$, there is $\epsilon'>0$ such that
\be{d4-top}
P\pare{  l_n(\D_n(y_n, n))\ge x_n/2}\le \exp\pare{-n^{\epsilon'} \xi_n^2}.
\ee
A bound of the type $P(-X_n\ge x_n)\le \exp(-c \xi_n^2)$ now
follows from \reff{neg.17}, and \reff{neg-prop3.3} after
we choose $\lambda$ small enough in the last term of
the right hand side of \reff{neg.10}.

\subsubsection{Proof of \reff{neg.1-bis}}
We fix $A$ large constant, and take the subdivision $\{b_1,\dots,b_M\}$
of $[A,y_n[$ with $b_1=A$, $b_{i+1}=2b_i$, for $i=1,\dots,M-1$,
with $M$ of order
$\log(n)$. We will choose $q$ slightly larger than 2, to be in
the super-critical regime (when $d\ge 4$), and we define
\be{def-Gi}
\G_i=\acc{ |\D_n(b_i,b_{i+1})|<\frac{C_1 n}{b_{i+1}^q}}.
\ee
Finally, for $q>2$, choose $p_i=p 2^{-i(q-2)/2}$ where $p$ is such that
$\sum_i p_i=1$. Now,
\be{d4-LT5}
\begin{split}
P\big(\sum_i\sum_{z\in \D_n(b_i,b_{i+1})}&l_n(z)(1-\zeta_z(l_n(z)))\ge
x_n\big)\le P\pare{\cup_i \G_i^c}\\
&\qquad +\sum_iP\pare{
\sum_{z\in \D_n(b_i,b_{i+1})}\frac{l_n(z)}{b_{i+1}}
(1-\zeta_z(l_n(z)))\ge \frac{x_n}{b_{i+1}},\ \G_i} .
\end{split}
\ee
First, we deal with $P(\cup_i \G_i^c)$ in the right hand side
of \reff{d4-LT5}. Note that
\be{cond-Gi}
\cup_i \G_i^c\subset \acc{ \|\ind_{\D_n(A,y_n)}l_n\|^q_q\ge
\frac{C_1}{2^q} n}.
\ee
We choose $C_1=2^{q+1} \kappa(q,d)$, and use Lemma~\ref{lem-sup}
to obtain, for any $\epsilon'>0$,
\be{neg.14}
P\pare{\cup_i \G_i^c} \le \exp\pare{-\frac{n^{1/q_c(d)-\epsilon'}}
{y_n^{q/q_c(d)-1}}}.
\ee
We neglect $P(\cup \G_i^c)$ if $n^{1/q_c(d)-\epsilon'}\ge
y_n^{q/q_c(d)-1} \xi_n^2$.  Since $\log(\xi_n)\le (d/2-\epsilon)/(d+4)
\log(n)$, and we are interested in $q$ close to 2, we only need to
check that taking $q=2$, for any $\epsilon>0$, we can find 
$\epsilon'>0$ such that
\be{check-1}
\frac{1}{q_c(d)}-\frac{1}{2}\pare{ \frac{2}{q_c(d)}-1}-\epsilon'\ge
\pare{ 2-( \frac{2}{q_c(d)}-1)}\pare{\frac{d/2-\epsilon}{d+4}}
\Longleftrightarrow \frac{1}{2}-\epsilon'\ge \frac{1}{2}-
\frac{\epsilon}{d}.
\ee
Since \reff{check-1} holds, 
we can find $\delta>0$ small enough, and $q=2+\delta$ so that
$P(\cup \G_i^c)$ is negligible.

We fix a realization of the random walk and integrate
first with respect to charges. For the charges, we
use the gaussian bounds of Remark~\ref{rem-hoffding} which
states that $\tilde \Gamma(x,n)\le \bar \chi_1x^2$, where
$\bar \chi_1=\max(1,\chi_1)$. In other words, on the event
$\G_i=\{|\D_n(b_i,b_{i+1})|\le C_1n/b_{i+1}^q\}$, we use
\be{d4-LT1}
\begin{split}
Q&\pare{\sum_{i=1}^M\sum_{z\in \D_n(b_i,b_{i+1})} l_n(z)
\pare{1-\zeta_z(l_n(z))}>\sum_i p_i x_n}\\
&\qquad\qquad\qquad\le 
\sum_{i=1}^M Q\pare{\sum_{z\in \D_n(b_i,b_{i+1})}
\frac{l_n(z)}{b_{i+1}}
\pare{1-\zeta_z(l_n(z))}> \frac{p_i}{b_{i+1}}x_n}.
\end{split}
\ee
Now, we consider a fixed $i\in \{1,\dots,M\}$, and on $\G_i$, we have for
any $\theta>0$
\be{d4-LT2}
\begin{split}
Q\pare{\sum_{z\in \D_n(b_i,b_{i+1})}\frac{l_n(z)}{b_{i+1}}
\pare{1-\zeta_z(l_n(z))}> \frac{p_i}{b_{i+1}}x_n}\le &
\exp\pare{-\frac{p_ix_n\theta }{
b_{i+1}}+\bar\chi_1 | \D_n(b_i,b_{i+1})|\theta^2}\\
\le &\exp\pare{-\frac{p_ix_n\theta }{b_{i+1}}+
\bar\chi_1 C_1\frac{n}{b_{i+1}^q}\theta^2}.
\end{split}
\ee
Note that if $|\D_n(b_i,b_{i+1})|\le p_i x_n/b_{i+1}$, then
the left hand side of \reff{d4-LT2} vanishes. Therefore, we 
assume that $|\D_n(b_i,b_{i+1})|> p_i x_n/b_{i+1}$, so that
the $\theta$ which minimizes the right hand side of \reff{d4-LT2}
is lower than 1, and we obtain
\be{d4-LT3}
P\pare{\sum_{z\in \D_n(b_i,b_{i+1})}\frac{l_n(z)}{b_{i+1}}
\pare{1-\zeta_z(l_n(z))}> \frac{p_i}{b_{i+1}}x_n,\ \G_i}\le 
\exp\pare{ -\frac{p_i^2 b_{i+1}^{q-2}\xi_n^2}{4C_1}}.
\ee
With our choice of $p_i,b_i$, we have that
$p_i^2 b_{i+1}^{q-2}\ge p^2 A^{q-2}$. Combining \reff{d4-LT1}
and \reff{d4-LT3}, we have
\be{d4-LT4}
P(\sum_{z\in \Z^d} l_n(z)1-\zeta_z(l_n(z))\ge x_n/2)\le M
\exp\pare{-\frac{p^2 A^{q-2}\xi_n^2}{4C_1}}.
\ee
The bound \reff{neg.1-bis} follows from \reff{d4-top} and
\reff{d4-LT4}.
\subsubsection{Dimension $d\ge 4$, and $\frac{d+2}{d+4}<\beta<1$.}
This corresponds to Region III of \cite{AC05}. 
We set $x_n=\xi_n$, $\zeta_n =\xi_n^{\frac{d}{d+2}}$, and $y_n=
\xi_n/\zeta_n$. Instead of \reff{neg.10}, we use 
\be{neg.22}
\begin{split}
P\big(-X_n\ge\xi_n\big)&
\le \P_0\pare{ l_n(\D_n(y_n^{1+\epsilon},n))\ge \frac{\xi_n}{2}}\\
&+
\P_0\pare{\|\ind_{\D_n(1,y_n^{1+\epsilon})} l_n\|_2^2\ge y_n \xi_n}
+\exp\pare{-\zeta_n y_n^{-\epsilon}(\lambda \xi_2-\lambda^2\chi_1)}.
\end{split}
\ee
Proposition 3.3 of \cite{AC05} yields that there is $\epsilon'>0$
such that
\be{neg.21}
\P_0\pare{l_n(\D_n(y_n^{1+\epsilon},n))\ge \frac{\xi_n}{2}}\le
\exp(-\xi_n^{{\frac{d}{d+2}}-\epsilon'}).
\ee
Now $\zeta_n^{\frac{d+4}{d+2}}\ge n$, and by Lemma~\ref{lem-sup},
for any $\epsilon$
\be{neg.20}
\P_0\pare{ \sum_{z\in \D_n(1,y_n^{1+\epsilon})} l_n^2(z)\ge 
\xi_n^{\frac{d+4}{d+2}}} \le 
\exp(-\frac{\xi_n^{(\frac{d+4}{d+2})(\frac{1}{q_c(d)}-\epsilon)}}
{y_n{(\frac{2}{q_c(d)}-1)}}).
\ee
The upper bound in \reff{neg.2} follows from \reff{neg.22},
\reff{neg.20}, and \reff{neg.21}.

\section{Lower Bounds.}\label{sec-lblt}
In realizing the lower bounds for 
Theorems \ref{prop-neg3}, \ref{prop-neg4-I}, and \ref{prop-neg4-II}, 
two strategies of the walk are distinguished: (i)
the walk is localized a time $T_n$ into a ball of radius $r_n$ with 
$r_n^2\ll T_n$, (ii) the walk roams freely. 
\subsection{On localizing the walk}\label{sec-local}
We introduce two sequences $\{T_n,r_n,\ n\in \N\}$.
We force the random walk to spend a time
$T_n$ in the ball centered at 0, of radius $r_n$, that we denote $B(r_n)$.

If $\tau_n=\inf\{n\ge 0: S(n)\not\in B(r_n)\}$, it is 
well known that for some constant $c_0$
\be{LT-lowkey}
\P_0(\tau_n>T_n)\ge \exp\pare{-c_0\frac{T_n}{|B(r_n)|^{2/d}}}.
\ee
Once the walk is forced to stay inside $B(r_n)$, we turn
to estimating the cost of $\{X_n<-x_n\}$. We then choose 
$\{T_n,r_n\}$ so as to match the cost with \reff{LT-lowkey}.

First, we need some relation between being localized
a time $T_n$ in a ball $B(r_n)$, and
visiting enough sites of $B(r_n)$ a time of order $T_n/|B(r_n)|$.
We have shown in \cite{A06} Proposition 1.4,
that in $d=3$,  for sequences $\{r_n,T_n\}$ going to
infinity with $r_n^d\le K T_n$, for some constant $K$,
there are positive constants $\delta_0$ and $\epsilon_0$,
independent of $r_n,T_n$ such that, for $n$ large enough
\be{prop1.4}
\P_0\pare{|\{z:\ l_{T_n}(z)>\delta_0 \frac{T_n}{|B(r_n)|}\}|
\ge \epsilon_0 |B(r_n)|}\ge
\frac{1}{2}\P_0\pare{\tau_n>T_n}.
\ee
Let $\RR_n$ be the set of sites visited by the random
walk before time $n$.
The only fact used in proving \reff{prop1.4} is an asymptotical bound
on $\P_0(|\RR_n|<n/\xi)$ for a fixed large $\xi$ and $n$ going to infinity.
Now, there is an obvious relation between $|\RR_n|$ and $\|l_n\|_q$
which reads as follows. For $q>1$
\be{link-Rln}
\pare{\frac{n}{|\RR_n|}}^{q-1}\le \frac{\|l_n\|_q^q}{n}.
\ee
Thus, from \reff{link-Rln} and \cite{A08} Theorem 1.1, we have
for $\xi^{q-1}>\kappa(q,d)$, and $q<q_c(d)$
\be{LD-range}
\P_0\pare{ |\RR_n|<\frac{n}{\xi}}
\le \P_0\pare{\|l_n\|_q^q\ge \xi^{q-1} n}\le
\exp(-c_1^+ \xi^{\frac{2}{d}}n^{1-\frac{2}{d}}).
\ee
Since $q_c(d)=\frac{d}{d-2}>1$, as soon as $d\ge 3$, \reff{LD-range}
is sufficient to obtain \reff{prop1.4} by following the proof of
\cite{A06}, and we omit the details. We now 
focus on the following set of sites
\be{def-G}
\G_n=\acc{z:\ \delta_0 \frac{T_n}{|B(r_n)|}\le l_{T_n}(z)\le
\frac{2 T_n}{\epsilon_0 |B(r_n)|}}.
\ee
Note that 
\[
|\{z:\ l_{T_n}(z)>\frac{2 T_n}
{\epsilon_0 |B(r_n)|}\}|\le \frac{\epsilon_0}{2} |B(r_n)|,
\]
so that $\{l_{T_n}>\delta_0 T_n/|B(r_n)|\}=\G_n\cup 
\{l_{T_n}>2 T_n/(\epsilon_0 |B(r_n)|)\}$, and
\be{LT-low8}
\P_0\pare{ |\G_n|\ge \frac{\epsilon_0}{2} |B(r_n)|}\ge
\P_0\pare{|\{z:\ l_{T_n}(z)>\delta_0 \frac{T_n}
{|B(r_n)|}\}|\ge \epsilon_0 |B(r_n)|}.
\ee
Now, in the scenario we are adopting, it will be easy to estimate
the contribution of sites of $\G_n$, which is a random set.
To use the notations of Proposition~\ref{prop-snagaev}, we define
for $z\in \Z^d$, $Y_z^{(n)}=l_n(z)(1-\zeta_z(l_n(z)))$.
We have, for $\delta>0$ small
\be{LT-low1}
\acc{\sum_{z\in \Z^d} Y_z^{(n)} \ge x_n}
\supset\acc{\sum_{z\in \G_n} Y_z^{(n)} \ge (1+\delta)x_n}\cap
\acc{\sum_{z\not\in \G_n} Y_z^{(n)} \ge-\delta x_n}.
\ee
When we integrate \reff{LT-low1} over the charges, we use that charges over 
disjoint regions are independent. Thus, we fix a realization of the
walk, and
\be{LT-ind}
Q\pare{\sum_{z\in \Z^d} Y_z^{(n)} \ge x_n}\ge
Q\pare{\sum_{z\in \G_n} Y_z^{(n)} \ge (1+\delta)x_n}
Q\pare{\sum_{z\not\in \G_n} Y_z^{(n)} \ge-\delta x_n}.
\ee
We first deal with the charges in $\G_n^c$.
We show using \reff{LB.10} that on 
$\B_n=\{\| l_n\|_2\le x_n n^{-\epsilon'}\}$, for $\epsilon'$
small, then
\be{notinG}
\begin{split}
\ind_{\B_n} Q\pare{ \sum_{z\not\in \G_n} Y_z^{(n)} 
\le-\delta x_n}&\le \ind_{\B_n} \frac{
\sum_{z\in \Z^d} E[(Y_z^{(n)})^2]}{(\delta x_n)^2}\\
&\le  \ind_{\B_n} \frac{\chi_1 \sum_{z\in \Z^d} l_n^2(z)}
{(\delta x_n)^2}\le 
\ind_{\B_n} \frac{\chi_1}{\delta^2n^{2\epsilon'}}.
\end{split}
\ee
Thus, from \reff{notinG} with $n$ large, we have
\be{notinG2}
\ind_{\B_n}
Q\pare{\sum_{z\not\in \G_n} Y_z^{(n)} \ge -\delta x_n} \ge
\frac{\ind_{\B_n}}{2}
\ee
From \reff{LT-low1} and \reff{notinG2}, we obtain, when integrating
only over the charges
\be{LT-low2}
\ind_{\B_n} Q\pare{ \sum_{z\in \Z^d} Y_z^{(n)} \ge x_n} \ge
\frac{\ind_{\B_n}}{2} Q\pare{
\sum_{z\in \G_n} Y_z^{(n)} \ge (1+\delta)x_n}.
\ee
Thus, after integrating over the walk
\be{LT-low3}
\begin{split}
2P\pare{ \sum_{z\in \Z^d} Y_z^{(n)}\ge x_n}+\P_0\pare{\B_n^c}&\ge
P\pare{ \sum_{z\in \G_n} Y_z^{(n)} \ge (1+\delta)x_n}\\
&\ge P\pare{|\G_n|\ge \frac{\epsilon_0}{2} |B(r_n)|,
 \sum_{z\in \G_n} Y_z^{(n)} \ge (1+\delta)x_n} .
\end{split}
\ee
Assume for a moment that $\P_0(\B_n^c)$ were negligible.
When integrating only over charges the last term of \reff{LT-low3}, 
we invoke Nagaev's Proposition~\ref{prop-snagaev},
applied to $\{Y_z^{(n)},\ z\in \G_n\}$. 
To simplify notations, we
assume henceforth that $T_n=n$ (though 
we can force the transient walk never to
return to $\G_n$ after time $T_n$, so that for $z\in \G_n$ we would
have $l_n(z)=l_{T_n}(z)$). Now, when we fix a realization
of the walk, we have easily from the equality \reff{LB.10}, for constants
$\chi_1$ and $\chi_4$
\be{LT-low0}
\chi_1 l_n^2(z)\ge E_Q[(Y_z^{(n)})^2]\ge 2 (l_n^2(z)-l_n(z))\quad\text{and}
\quad E_Q[(Y_z^{(n)})^4]\le  \chi_4 l_n^4(z). 
\ee
From Jensen's inequality, we have $E_Q[|Y_z^{(n)}|^3]\le
\chi_3 l_n^3(z)$ with $\xi_3=\xi_4^{3/4}$. 
Note that in order to have a non-zero lower bound for the variance
of $Y_z^{(n)}$, we impose
\be{cond-appli3}
\delta_0 \frac{T_n}{|B(r_n)|}\ge 2\quad\text{so that}\quad\forall z
\in \G_n\quad E_Q[Y_z^2]\ge 2(l_n^2(z)-l_n(z))\ge l_n^2(z).
\ee

With the notations 
of Proposition~\ref{prop-snagaev}, we have (using \reff{LT-low0})
on $\{|\G_n|\ge \frac{\epsilon_0}{2} |B(r_n)|\}$
\be{LT-low4}
\frac{\epsilon_0 \delta_0^2}{2} \frac{T_n^2}{|B(r_n)|}\le \sigma_n^2\le
\frac{4\chi_1}{\delta_0\epsilon_0^2}\frac{T_n^2}{|B(r_n)|}\quad
\text{and}\quad \C_n^3\le \frac{8\chi_3}{\delta_0\epsilon_0^3}
\frac{T_n^3}{|B(r_n)|^2}.
\ee
Also, $\sigma_n t_n=(1+\delta)x_n$, so that \reff{cond-nagaev}
holds if for some $\epsilon_\NN>0$, and $n$ large enough
\be{cond-appli1}
\sigma_n\le (1+\delta)x_n,\quad
(1+\delta)x_n \C_n^3\le \epsilon_\NN \sigma_n^4,\quad
\text{and}\quad
(1+\delta)x_n \max_{z\in \G_n} \sqrt{ E\cro{(Y_z^{(n)})^2}}
\le \epsilon_\NN \sigma_n^2.
\ee
Using \reff{LT-low4}, \reff{cond-appli1} and \reff{cond-appli3}
follow if, for some constant $c_1$ 
\be{cond-appli2}
\frac{4\chi_1}{\delta_0\epsilon_0^2}\frac{T_n^2}{|B(r_n)|}\le x_n^2,
\quad\text{and}\quad 
x_n\le \epsilon_\NN c_1 T_n.
\ee
When \reff{cond-appli2} holds, and we can use 
Proposition~\ref{prop-snagaev},
to obtain on $\{|\G_n|\ge \frac{\epsilon_0}{2} |B(r_n)|\}$,
and for constants $c_1,c_2$
\be{snagaev-appli1}
Q\pare{ \sum_{z\in \G_n} Y_z^{(n)} \ge (1+\delta)x_n}
\ge \exp\pare{-c_1 (\frac{x_n}{\sigma_n})^2}\ge
\exp\pare{-c_2 \frac{x_n^2\ |B(r_n)|}{T_n^2}}.
\ee
After integrating over the walk, recalling \reff{prop1.4},
\reff{LT-lowkey}, \reff{LT-low3} and \reff{LT-low8}, we have
\be{LT-low7}
\begin{split}
2P\pare{ \sum_{z\in \Z^d} Y_z^{(n)}\ge x_n}
&\ge P\pare{|\G_n|\ge \frac{\epsilon_0}{2} |B(r_n)|,   
\sum_{z\in \G_n} Y_z^{(n)} \ge (1+\delta)x_n}-\P_0(\B_n^c)\\
&\ge \exp\pare{ -c_2 \frac{x_n^2\ |B(r_n)|}{T_n^2}-c_0\frac{
T_n}{|B(r_n)|^{2/d}}}-\P_0(||l_n||_2^2\ge x_n^2 n^{-2\epsilon'}).
\end{split}
\ee
From inequality \reff{LT-low7}, the difference between $d=3$ and
$d\ge 4$ is obvious, when imposing a localisation of the walk.
Indeed, matching the two costs in \reff{LT-low7}, we find
\be{2-costs}
\frac{x_n^2 |B(r_n)|}{T_n^2}=\frac{T_n}{|B(r_n)|^{2/d}}
\Longrightarrow |B(r_n)|^{\frac{d+2}{d}}=\frac{T_n^3}{x_n^2}.
\ee
Thus, combining \reff{LT-low7} with the choice of \reff{2-costs},
we obtain for a constant $c_d^->0$
\be{lower-cost}
P(X_n\le -x_n)\ge \exp\pare{-c_d^-x_n^{\frac{4}{d+2}}T_n^{\frac{d-4}{d+2}}}
-\P_0(\B_n^c).
\ee
Corollary~\ref{cor-2beta} shows that $\P_0(\B_n^c)\ll 
\exp(-c_d^- \xi_n^{\frac{d}{d+2}})$. Henceforth, we neglect
$\P_0(\B_n^c)$. 

\subsection{The case $d=3$ and $a_0\le \xi_n\le n^{1/3}$.}
\label{sec-LBd3}
In this section, we choose $T_n=n$, and $|B(r_n)|^{5/3}=n^{3}/x_n^{2}$,
as suggested in \reff{2-costs}.

We start with $\xi_n\le c_1\epsilon_\NN  n^{1/3}$.
In this case, $x_n=\xi_n n^{2/3}$. The discussion of the previous
section applies here.
Note that sites of $\G_n$ are visited about $\xi_n^{6/5}$-times each.
Conditions~\reff{cond-appli2} are satisfied, and the discussion following
it holds. The bound \reff{lower-cost} provides the desired lower bound.

Now, we deal with $x_n=\xi n$, with $1>\xi\ge c_1\epsilon_\NN$.
The second inequality in \reff{cond-appli2} fails, and Nagaev's
lower bound cannot be applied.
We choose $\delta>0$ small enough so that $\xi(1+\delta)^2<1$, and
we consider the event $\A=\{\forall z\in B(r_n), (1-\zeta_z)\ge
\xi(1+\delta)^2\}\cap\{\tau_n>n\}$. Note that
\[
\A\subset \acc{\sum_{z\in \Z^d} l_n(z)(1-\zeta_z(l_n(z))
\ge \xi(1+\delta)^2 n}.
\]
However, there might be some sites of $B(r_n)$ that
the walk visits once, and if $\eta\in \{-1,1\}$, we will have
on this sites that $\zeta_z(l_n(z))=0$. We will restrict
to sites of $B(r_n)$ visited often. Note that, for $\alpha(\xi)>0$,
\[
\lim_{n\to\infty} Q\pare{ 1-\zeta_z(n)\ge \xi(1+\delta)^2}=
\lim_{n\to\infty} Q\pare{ \pare{\frac{1}{\sqrt n} \sum_{i=1}^n
\eta_z(i)}^2\le 1-\xi(1+\delta)^2}=\alpha(\xi).
\]
Thus, there is $n_1$ (depending on $\xi$ and $\delta$)
such that for $n\ge n_1$
\[
Q\pare{ \pare{\frac{1}{\sqrt n} \sum_{i=1}^n
\eta_z(i)}^2\le 1-\xi(1+\delta)}\ge \frac{1}{2} \alpha(\xi).
\]
Now, with $n_1$ fixed, we define a set
\[
\G_n=\acc{z\in B(r_n):\ l_n(z)\ge n_1}.
\]
On the event $\{\tau_n>n\}$, we have for $n$ large enough
(using that $|B(r_n)|\ll n$)
\[
l_n(\G_n^c)\le |B(r_n)| n_1\Longrightarrow
l_n(\G_n)\ge n-|B(r_n)| n_1\ge \frac{n}{1+\delta}.
\]
Thus,
\[
\A\subset \acc{\sum_{z\in \G_n} l_n(z)(1-\zeta_z(l_n(z))\ge 
l_n(\G_n) \xi(1+\delta)^2=\xi(1+\delta) n}
\]
Using \reff{LT-low3} (with $\delta$ occurring in \reff{LT-low3}),
 we have
\be{LT-33}
\begin{split}
2P\pare{\sum_{z\in \Z^d} Y_z\ge \xi n}+\P_0\pare{\B_n^c}\ge&
(\frac{\alpha(\xi)}{2})^{|B(r_n)|}\times\P_0\pare{\tau_n>n}\\
\ge&(\frac{\alpha(\xi)}{2})^{|B(r_n)|}\times
\exp\pare{-c_0 \frac{n}{|B(r_n)|^{2/d}}}.
\end{split}
\ee
Since $1>\xi>c_1 \epsilon_\NN$, the power of $\xi$ appearing in
\reff{LT-33} is irrelevant. We only need to check
that the speed exponent is correct in \reff{LT-33}

\subsection{The case $d\ge 4$ and $n^{\frac{d+2}{d+4}}\ll \xi_n\ll n$}
\label{sec-scen2}
Here $x_n=\xi_n$.
Assume that we localize the walk a time $T_n$ inside $B(r_n)$.
We make use of Section~\ref{sec-local} until the point
where we assumed $T_n=n$ (that is a paragraph before \reff{LT-low0}).
If we were allowed to
identify the two costs in \reff{LT-low7}, we would find here
$T_n=x_n=\xi_n$, and $|B(r_n)|=\xi_n^{\zeta_d}$,
 with $\zeta_d=\frac{d}{d+2}$.
Note that in dimension 4 or larger, with $T_n$ of order $\xi_n$,
we are not entitled to use Nagaev's lower bound. 
On the other hand, $|B(r_n)|=\xi_n^{\zeta_d}$, is the expected speed,
so that
constraining the local charges on $\G_n$ would yield the correct
cost. We observe that
we are entitled to use the CLT for $\zeta_z(l_n(z))$,
for each sites in $\G_n$, since $l_n(z)\ge l_{T_n}(z)\ge\xi_n^{1-\zeta_d}$.
With the notation $Z$ for a 
standard gaussian variable, and $n$ 
large enough, we have for $z\in \G_n$,
and uniformely over $l_n(z)$
\[
\alpha_0:=\frac{1}{2}P(Z^2<\frac{1}{2})\le
Q(\zeta_z(l_n(z))< \frac{1}{2}).
\]
With the choice $T_n=\frac{4}{\epsilon_0}\xi_n$
(note that $T_n\ll n$ for $n$ large),
recalling the definition of $\G_n$ in \reff{def-G}, and
using that $l_n(z)\ge l_{T_n}(z)$
\[
\acc{\forall z\in \G_n, \zeta_z(l_n(z))<
\frac{1}{2}}\cap \acc{|\G_n|\ge \frac{\epsilon_0}{2}|B(r_n)|}
\subset \acc{\sum_{z\in \G_n} Y_z\ge \frac{1}{2}|\G_n|T_n=
(1+\delta) \xi_n}.
\]
Thus, using \reff{LT-low3}
\be{LT-low10}
2P\pare{\sum_{z\in \Z^d} Y_z\ge \xi_n}+\P_0\pare{\B_n^c}\ge
\alpha_0^{|B(r_n)|}\times\P_0\pare{\tau_n>T_n}
\ge \exp\pare{-c_d^- \xi_n^{\zeta_d}}.
\ee
\subsection{The case $d\ge 4$ and $x_n=\xi n$}
We assume that $\xi<1$, for $\delta'>0$ so small that
$(1+\delta')\xi<1$, we choose $T_n=(1+\delta')\xi n$ and 
$|B(r_n)|=(\xi n)^{d/(d+2)}$. 
We force the local charges to realize $1-\zeta_z(l_n(z))\ge
1-\frac{\delta'}{4}$ for $\delta'$ arbitrarily small. Note that for $\alpha_1>0$,
\[
\lim_{n\to\infty} Q\pare{ \pare{\frac{1}{\sqrt n} \sum_{i=1}^n
\eta_z(i)}^2\le \frac{\delta'}{4}}=\alpha_1.
\]
Thus, there is $n_1$ (depending on $\xi$ and $\delta'$)
such that for $n\ge n_1$
\be{LT-low5}
Q\pare{ \pare{\frac{1}{\sqrt n} \sum_{i=1}^n
\eta_z(i)}^2\le \frac{\delta'}{4}}\ge \frac{1}{2} \alpha_1.
\ee
Now, using $n_1$, we define a set
\[
\G_n=\acc{z\in B(r_n):\ l_n(z)\ge n_1}.
\]
On the event $\{\tau_n>(1+\delta')\xi n\}$, we have for $n$ large enough
(using that $|B(r_n)|\ll n$)
\[
l_n(\G_n^c)\le |B(r_n)| n_1\Longrightarrow
l_n(\G_n)\ge (1+\delta')\xi n-|B(r_n)| n_1\ge (1+\delta')
(1-\frac{\delta'}{4})n\xi. 
\]
We use \reff{LT-low5} for $\zeta_z(l_n(z))$, with $z\in \G_n$.
Thus, on $\{\tau_n\ge (1+\delta')\xi n\}$,
\[
\acc{\forall z\in \G_n, \zeta_z(l_n(z))<\delta'}
\subset \acc{\sum_{z\in \G_n} Y_z\ge (1-\frac{\delta'}{4}) 
l_n(\G_n)\ge (1+\delta')(1-\frac{\delta'}{4})^2n\xi}.
\]
Now, we choose $\delta'$ so small 
that $(1+\delta')(1-\frac{\delta'}{4})^2
\ge 1+\delta$, for $\delta$ occurring in \reff{LT-low3}.
Thus, using \reff{LT-low3}
\be{LT-34}
\begin{split}
2P\pare{\sum_{z\in \Z^d} Y_z\ge \xi_n}+\P_0\pare{\B_n^c}&\ge
(\frac{\alpha_1}{2})^{|B(r_n)|}\times\P_0\pare{\tau_n>T_n}\\
& \ge (\frac{\alpha_1}{2})^{|B(r_n)|}\times 
\exp\pare{-c_0 \frac{(1+\delta')\xi n}{|B(r_n)|^{2/d}}}\\
&\ge \exp(-c_d^-(\xi n)^{\frac{d}{d+2}}).
\end{split}
\ee
This yields the desired bound.
\subsection{The case $d\ge 4$ and $n^{2/3}\ll \xi_n\ll n^{(d+2)/(d+4)}$.}
The strategy in this region
(region I of \cite{AC05}) consists in letting the walk roam
freely, while the {\it local charges} perform a moderate deviations.
Note that our scenery $\zeta_z$ depends on the local times, and
on sites visited only once by the walk, $Y_z$ may vanish
by \reff{LB.10}, as in the model where $\eta\in \{-1,1\}$.
Thus, we only consider sites where $\{z: l_n(z)=2\}$, since
$\frac{1}{2}(\eta_1+\eta_2)^2-1$ is not degenerate. Also, 
a transient random walk has enough sites of this type. Indeed,
Becker and K\"onig in \cite{BK} have shown that, in $d\ge 3$ 
with $\D_n(k)=\{z: l_n(z)=k\}$ for integer $k$, we have
\be{limit-BK}
\lim_{n\to\infty} \frac{E\cro{|\D_n(k)|}}{n}=\gamma_0^2(1-\gamma_0)^{k-1},
\quad\text{where}\quad \gamma_0=\P_0(S(k)\not= 0,\ \forall k>0).
\ee
We choose a scenario based only on $\D_n(2)$. Note that 
for $n$ large enough,
the fact that $|\D_n(2)|\le n$, and \reff{limit-BK} imply that 
\[
\frac{1}{2} \gamma_0^2(1-\gamma_0)\le \frac{E\cro{|\D_n(2)|}}{n}\le
\P_0\pare{ \frac{|\D_n(2)|}{n}\ge \frac{1}{4} \gamma_0^2(1-\gamma_0)}
+\frac{1}{4} \gamma_0^2(1-\gamma_0).
\]
Thus,
\be{D2-big}
\P_0\pare{ \frac{|\D_n(2)|}{n}\ge \gamma_1}
\ge \gamma_1\quad\text{with}\quad \gamma_1=
\frac{1}{4} \gamma_0^2(1-\gamma_0).
\ee
Now, we consider the following decomposition, for $\delta>0$ small
(recall that here $x_n=\sqrt{n}\ \xi_n$)
\be{LT-low11}
\acc{\sum_{z\in \Z^d} Y_z^{(n)} \ge \sqrt{n}\ \xi_n}
\supset\acc{\sum_{z\in \D_n(2)} Y_z^{(n)} \ge (1+\delta)\sqrt{n}\ \xi_n}\cap
\acc{\sum_{z\not\in \D_n(2)} Y_z^{(n)} \ge-\delta\sqrt{n}\ \xi_n}.
\ee
We treat the second event on the right hand side of \reff{LT-low11}
as in Section~\ref{sec-local}: we restrict to $\B_n$ 
(where $P(\B_n^c)$
is negligible by Corollary~\ref{cor-2beta}), and we use Markov's inequality.

Now, fixing a realization of the walk,
$\{Y_z,\ z\in \D_n(2)\}$ are centered i.i.d with $E[Y_z^2]=2(
E_Q[\eta^4]+1)$, and on $\{|\D_n(2)|> \gamma_1
n\}$, then $\{\sum_{\D_n(2)} Y_z\ge(1+\delta)\sqrt{n}\ \xi_n\}$ is 
a moderate
deviations. Thus, there is a constant $\sous c$, 
such that on the event $\{|\D_n(2)|> \gamma_1n\}$, and for $n$ large
\be{LT-low12}
\begin{split}
Q\pare{ \sum_{\D_n(2)} Y_z\ge (1+\delta)\sqrt{n}\ \xi_n}&\ge
\sous c 
\exp\pare{-\frac{((1+\delta)\xi_n)^2n}{2|\D_n(2)|(E_Q[\eta^4]+1)}}\\
&\ge \sous c \exp\pare{-\frac{(1+\delta)^2\xi_n^2}
{2\gamma_1(E_Q[\eta^4]+1)}}.
\end{split}\ee
After integrating \reff{LT-low12} the walk's law, we have
\be{LT-low13}
P\pare{|\D_n(2)|> \gamma_1n,\ 
\sum_{\D_n(2)} Y_z\ge (1+\delta)\sqrt{n}\ \xi_n}\ge
\sous c  \gamma_1 \exp\pare{-\frac{(1+\delta)^2}
{2\gamma_1(E_Q[\eta^4]+1)}\xi_n^2}.
\ee
\section{Proof of Proposition~\ref{prop-gibbs}}\label{sec-gibbs}
\paragraph{Large $\beta$}
First, $H_n\ge -n$ implies the upper bound in \reff{gibbs-1}.
The lower bound in \reff{gibbs-1} follows from the lower bound
in \reff{neg.3} with $\xi_n=\xi n^{1/3}$, 
and the following inequalities: for $\xi<1$
\be{gibbs-7}
\begin{split}
Z_n^-\pare{\frac{\beta}{n^{2/5}}}=&E\cro{\exp\pare{-\beta
\frac{H_n}{n^{2/5}}}}\ge P(H_n\le -\xi n) e^{\beta \xi n^{3/5}}\\
\ge & \exp\pare{n^{3/5}(\beta \xi -c^-_3\xi^{4/5})}.
\end{split}
\ee
For any fixed $\xi<1$, we choose $\beta$ large enough so that 
the lower bound in \reff{gibbs-1} holds.

Now, define
\[
\A_n(a)=\acc{|\{z\in \Z^d: \frac{ n^{\frac{2}{5}}}{a}\le l_n(z)\le
 a n^{\frac{2}{5}}\}|\ge  \frac{n^{3/5}}{a^4}}.
\]
Using the estimates of Section~\ref{sec-forgibbs}, we have for
$\chi>0$
\be{gibbs-6}
E\cro{\exp(\pare{-\beta\frac{H_n}{n^{2/5}}}}
\le e^{\beta n^{3/5}}P(\A_n^c(a))\le e^{n^{3/5}(\beta -\chi a^{2/3})}.
\ee
Choosing $a$ large enough so that $2\beta<\chi a^{2/3}$, and
using the lower bound in \reff{gibbs-7}, we obtain \reff{gibbs-2}.
\paragraph{Small $\beta$.}
First, we decompose the partition function over the three regimes
for $-H_n$: the moderate deviation, the large deviation and
intermediate regimes. Thus,
\be{gibbs-20}
Z_n^-(\frac{\beta}{n^{2/5}})=Z_I(\beta)+Z_{I\!I}(\beta)+
Z_{I\!I\!I}(\beta),
\ee
with for $\epsilon$ small
\[
Z_I(\beta)=E\cro{\exp\pare{-\beta\frac{H_n}{n^{2/5}}}\ind\acc{
n^{\frac{1}{2}+\epsilon}<-H_n<n^{\frac{2}{3}+\epsilon}}},
\]
\[
Z_{I\!I}(\beta)=E\cro{\exp\pare{-\beta\frac{H_n}{n^{2/5}}}\ind\acc{
n^{\frac{2}{3}+\epsilon}<-H_n<n}},
\]
and $Z_{I\!I\!I}(\beta)$ correponds to the remaining regimes.

We first deal with $Z_I(\beta)$ and rely on Chen's result \reff{intro.4}.
We note that from Chen's proof,
his asymptotic result of \reff{intro.4} is actually
uniform in the sequence $\xi_n$, in the sense that there is a
sequence $\{\delta_n\}$ going to 0, such that for any
$\xi_n\in [n^\epsilon,n^{1/6-\epsilon}]$, we have
\be{chen-turbo}
P(\frac{-H_n}{\sqrt n}>\xi_n)=\exp\pare{-\frac{\xi_n^2}{2c_d}(1+\delta_n)}.
\ee
We have
\be{gibbs-23}
\begin{split}
Z_I(\beta)=& \exp( \beta n^{1/10+\epsilon})+
\beta\int_{n^{1/10+\epsilon}}^{n^{4/15-\epsilon}}e^{\beta u}
P\pare{\frac{-H_n}{n^{2/5}}>u}\ du\\
=& \exp( \beta n^{1/10+\epsilon})+ \beta n^{1/10}
\int_{n^\epsilon}^{n^{1/6-\epsilon}}
\exp\pare{ \beta n^{1/10} u-\frac{u^2}{2c_d}(1+\delta_n)}\ du
\end{split}
\ee
Now, the aymptotic behaviour is found as we maximize
$\beta n^{1/10} u-\frac{u^2}{2c_d}$, which is $c_d\beta^2 n^{1/5}/2$.
In other words, it is a simple computation that we omit, which
yields for any $\beta>0$,
\be{gibbs-24}
\lim_{n\to\infty} \frac{1}{n^{1/5}}\log Z_I(\beta)=
\ \frac{ c_d \beta^2}{2}.
\ee
We deal now with $Z_{I\!I}$, which corresponds to regime studied in
Theorem~\ref{prop-neg3}. We will show that for $\beta$ small,
$Z_{I\!I}(\beta)\le \exp( \epsilon n^{1/5})$, for $\epsilon$ small.
Note that
\be{gibbs-25}
Z_{I\!I}(\beta)\le \sum_{k=0}^{\log_2(n^{1/3})}
e^{2^{k+1}n^{4/15+\epsilon}\beta}\ 
P\pare{ 2^{k}n^{4/15+\epsilon}\le \frac{-H_n}{n^{2/5}}<
2^{k+1}n^{4/15+\epsilon}}
\ee
In view of \reff{gibbs-25},
it is enough to show that for $n^{3/5}\ge \xi_n\ge n^{4/15+\epsilon}$,
we have
\be{gibbs-8}
P(-H_n\ge \xi_n n^{2/5})\le e^{-2 \beta\xi_n}.
\ee
From \reff{neg.3}, we have in this regime
\be{gibbs-9}
P(-H_n\ge \xi_n n^{2/5})\le \exp\pare{ -c_3^+
\pare{ \xi_n n^{2/5-2/3}}^{4/5}n^{1/3}},
\ee
and \reff{gibbs-8} requires that 
\be{gibbs-10}
c_3^+ \xi_n^{4/5} n^{3/25}\ge 2\beta \xi_n\Longleftrightarrow
\xi_n\le \pare{\frac{c_3^+}{2\beta}}^5 n^{3/5}.
\ee
Since $\xi_n\le n^{3/5}$, \reff{gibbs-10} holds if $\beta<c_3^+/2$.

Finally, we deal with  $Z_{I\!I\!I}$.
\be{gibbs-11}
\begin{split}
Z_{I\!I\!I}\le& \exp\pare{ \beta n^{1/2-2/5+\epsilon}}+
\exp(\beta n^{2/3-2/5+\epsilon})P(-H_n\ge n^{2/3-\epsilon})\\
\le & \exp\pare{ \beta n^{1/10+\epsilon}}+
\exp\pare{ -\frac{n^{1/3-\epsilon}}{4c_d}+\beta n^{4/15+\epsilon}}.
\end{split}
\ee
$Z_{I\!I\!I}$ is negligible when $\epsilon$ is such that
$\frac{4}{15}+3\epsilon\le \frac{1}{3}$.

We finally show \reff{gibbs-4}.
We choose $p>1$ such that $p\beta<\beta_1$, and use H\"older's inequality
\be{gibbs-12}
\begin{split}
E\cro{ e^{-\beta\frac{H_n}{n^{2/5}}}\ind_{\acc{l_n(z)>b n^{1/5}}
\not=\emptyset}}\le &
\pare{E\cro{e^{-p\beta\frac{H_n}{n^{2/5}}}}}^{1/p}\pare{
P(\exists z,\ l_n(z)>b n^{1/5})}^{1/q}\quad (q=\frac{p}{p-1})\\
\le & e^{C\beta^2 n^{1/5}} \pare{ n P_0(l_n(0)>b n^{1/5})}^{1/q}\\
\le & n^{1/q} \exp\pare {(C\beta^2-\frac{\chi_db }{q})n^{1/5}}.
\end{split}
\ee
As we choose $b$ large enough in \reff{gibbs-12}, we obtain
\reff{gibbs-4}.

\end{document}